\documentclass[12pt,a4paper]{amsart}
\usepackage{url}
\sloppypar
\usepackage[ansinew]{inputenc}	%
\usepackage{graphics}
\usepackage{a4}
\usepackage{amssymb}
\usepackage{amsmath}
\usepackage{amstext}
\usepackage{amsopn}
\usepackage{amsthm}
\usepackage{graphics, graphicx}
\usepackage{algpseudocodex}
\newtheorem{thm}{Theorem}[section]

\newtheorem{lem}[thm]{Lemma}

\theoremstyle{definition}
\theoremstyle{definition}
 
\newtheorem{de}[thm]{Definition}
\newtheorem{alg}[thm]{Algorithm}

\newtheoremstyle{uebstyle}{}{}{\footnotesize}{}{\normalsize\scshape}{}{\newline}{}
\theoremstyle{uebstyle}
\newtheorem{uebungen}[thm]{Übungsaufgaben}
\newtheorem{loesungen}[thm]{Lösungen zu Übungsaufgaben}
\numberwithin{equation}{section}

\newcommand{\mycomment}[1]{}

\DeclareMathAlphabet{\mathsc}{OT1}{cmr}{m}{sc}
\DeclareMathAlphabet{\mathbfsl}{OT1}{cmr}{bx}{it}

\newcommand{\vb}[1]{\mathbfsl{#1}}

\newcommand{\VecTwo}[2]
           { \left( 
             \begin{smallmatrix}
                           {#1} \\ {#2} 
             \end{smallmatrix} 
            \right) }
\newcommand{\N}{\mathbb{ N}} 
\newcommand{\Z}{\mathbb{ Z}} 
\newcommand{\R}{\mathbb{ R}} 
\newcommand{\Q}{\mathbb{ Q}}

\newcommand{\setsuchthat}{\,\,\pmb{|}\,\,}

\newcommand{\ul}[1]{\underline{#1}}

\newcommand{\sgn}{\mathrm{sgn}}

\newcommand{\LT}{\mathsc{Lt}}
\newcommand{\lt}{\operatorname{lt}}
\newcommand{\LM}{\mathsc{Lm}}
\newcommand{\lm}{\operatorname{lm}}
\newcommand{\LC}{\mathsc{Lc}}

\newcommand{\DegD}{\md_{\delta}}
\newcommand{\Min}{\operatorname{Min}}
\newcommand{\lea}{\le}
\newcommand{\gea}{\ge}
\newcommand{\strictlea}{<}

\renewcommand{\div}{\sqsubseteq}
\newcommand{\divd}{\div_{\delta}}
\newcommand{\md}{\mathsc{Deg}}
\newcommand{\NON}{{\N_0}^n}
\newcommand{\NONk}{\NON \times \ul{k}}
\newcommand{\dhat}{\hat{\delta}}

\newcommand{\ups}[1]{#1{\uparrow}}

\newcommand{\kvx}{k[\vb{x}]}

\newcommand{\rxn}{R[x_1,\ldots,x_n]}
\newcommand{\rxnk}{R[x_1,\ldots,x_n]^k}
\newcommand{\rvx}{R[\vb{x}]}
\newcommand{\rvxk}{R[\vb{x}]^k}

\newcommand{\submod}[1]{\langle #1 \rangle}

\newcommand{\term}[3]{#1 \vb{x}^{#2} e_{#3}}
\newcommand{\mono}[1]{\vb{x}^{#1}}
\newcommand{\monov}[2]{\term{}{#1}{#2}}
\newcommand{\termp}[2]{#1 \vb{x}^{#2}}
\newcommand{\monoexp}[2]{(#1, #2)}
\newcommand{\monosnk}{\operatorname{Mon} (n,k)}

\newcommand{\col}{\operatorname{col}}
\newcommand{\row}{\operatorname{row}}

 \renewcommand{\emptyset}{\varnothing}
\date{\today}
\usepackage[colorlinks=true,linkcolor=blue,filecolor=blue,citecolor=blue,urlcolor=blue]{hyperref}				%

\begin{document}
\bibliographystyle{is-alpha}

\title[Strong Gr\"obner bases over Euclidean domains]
      {Strong Gr\"obner bases and linear algebra in
       multivariate polynomial rings
       over Euclidean domains}
\thanks{Supported by the Austrian Science Fund (FWF):P33878}
\begin{abstract}
  We provide a self-contained introduction to Gr\"obner bases
  of submodules of $\rxnk$, where $R$ is a Euclidean domain,
  and explain how to use these bases to solve linear systems
  over $\rxn$.
\end{abstract}
\author{Erhard Aichinger}
\address{Erhard Aichinger,
Institute for Algebra,
Johannes Kepler University Linz,
4040 Linz,
Austria}
\email{\tt erhard@algebra.uni-linz.ac.at}

\maketitle
\section{Introduction}
The computation of Gr\"obner bases is a broadly applicable
method that solves many questions involving polynomials.
One such question is solving systems of linear equations over
a commutative ring $D$. Here, for a matrix $A \in D^{r \times s}$ and
$b \in D^r$, one would like to compute an $x \in D^s$ with $Ax = b$
(when it exists) and a basis of the module
$\ker (A) = \{x \in D^s \mid Ax = 0\}$;
then $\{x+k \mid k \in \ker (A)\}$ is the set of solutions of the
linear system. When $D$ is a multivariate polynomial ring
such as $\Z[x_1, x_2]$ or $\Q[x_1, x_2]$, Gr\"obner bases are
a tool to solve these questions. Being able to solve linear
systems over $D$ allows us to determine ideal membership in $D$
(solve $d_1x_1 + \cdots + d_n x_n = d$ over $D$ to find out
whether $d$ lies in the ideal generated by $d_1, \ldots, d_n$)
and to compute least common multiples (solve $x_1 - d_1 x_2 = x_1 - d_2 x_3 = 0$
for finding $x_1$ as a common multiple of $d_1, d_2$), and hence
also greatest common divisors when $D$ is a unique factorization domain.

When $D$ is a field, solving linear equations is accomplished by Gau\ss's algorithm. When $D = \Z$, solving systems of linear diophantine equations
can be done computing the Hermite normal form of a matrix.
Both cases are contained in the case that $D = \rxn =: \rvx$, where
$R$ is a Euclidean domain, and in the present note, we explain how
to solve linear systems over $\rvx$. The role of
the row echelon form of a matrix (when $D$ is a field)
and of the Hermite normal form (when $D$ is a Euclidean
domain) will be taken by a matrix whose rows are a Gr\"obner basis
of the module generated by the rows of the matrix.
Our approach includes computing
Gr\"obner bases over $\Z$, allowing us to do linear
algebra over $\Z[x_1, \ldots, x_n]$. As every finitely generated ring
is isomorphic to a quotient $\Z[x_1, \ldots, x_n]$ by an
ideal $I$ (which is finitely generated by Hilbert's Basis Theorem),
this will allow us to solve linear systems over all finitely
generated rings. In fact, the Gr\"obner basis algorithm presented
here (which is a modification of the algorithm given in
\cite{Li:ECOS}) will contain both Gau\ss's algorithm and
the Hermite normal form as special instances.

As we do not presuppose any knowledge on Gr\"obner bases, let us
start with a rough description: Given a submodule
$M$ of the $D$-module $D^s$, a Gr\"obner basis is a set of generators
of $M$ with particularly useful properties.
They were introduced by B.\ Buchberger, who presented an algorithm to compute such
bases when $D = k[x_1, \ldots, x_n]$ ($k$ a field) and $s=1$ \cite{Bu:EAZA,Bu:EAKF}
and named them in honour of his supervisor W.\ Gr\"obner. Generalizing to $s>1$ is
then straightforward.
The case $D = \Z[x_1, \ldots, x_n]$ provides additional difficulties,
as now diophantine linear equations over $\Z$ are included.
For this case, several types of Gr\"obner bases were introduced.
We will use ``strong Gr\"obner bases'', and one main source for
our development is \cite{Li:ECOS}.
The case that $D = R[x_1, \ldots, x_n]$ for a Euclidean
domain $R$ has been been treated in \cite{KK:CAGB}.
What our treatment adds to these is that we:
\begin{itemize}
\item consider bases of submodules of $D^s$ also for the case $s > 1$,
\item provide self-contained proofs of
  the Gr\"obner basis criterion~Theorem~\ref{thm:Spols}
   and  
  of the uniqueness of reduced strong Gr\"obner
  bases,
\item explain how reduction and augmentation by $S$-polynomials
     can be interleaved in the course of the algorithm, and
\item explicitly state how to solve linear systems over $D$.
\end{itemize}
The computations will start from a generating
set $F$ of a submodule of $D^s$ and compute a set
of generators $G$ (a Gr\"obner basis) with certain desirable properties.
The algorithm is a sequence of
the following steps, which we illustrate here by examples
for the case $D = \Z[x,y]^3$.
\begin{itemize}
\item Augmentation: When $f = (10 x^2 y^2+y,0,x)$
  and $g = (4 x^3 y+x^2,1,0)$ are in $F$, then add
  $h = x \, f - 2 y \, g = (2 x^3 y^2-2 x^2 y+x y,-2 y,x^2)$ to $F$.
  One important property of $h$ is that the
  leading coefficient $2$ of $h$  is
  smaller than the leading coefficients $10$ of $f$ and $4$ of $g$.
  Such an $h$ is called an \emph{$S$-polynomial vector} (from ``subtraction'').
\item Reduction: When $f = (10 x^2 y^2+y,0,x)$ and
  $g = (x-2 y,1,0)$ are in $F$, then
  replace $f$ by $f' = f - 10 x y^2 g =  (20 x y^3+y,-10 x y^2,x)$.
\end{itemize}
When these simple steps are performed in some proper order, the process
will eventually
terminate and produce a set of generators $G$ that will have
the required properties.
Termination is proved using the fact that certain ordered sets
have no infinite descending chains.
The fact that the
final result $G$ has the desired properties uses a central theorem
on $S$-polynomials vectors. For the case that $R$ is a field
this theorem goes back to
\cite{Bu:EAZA} (cf. \cite[Theorem~3.3]{Bu:ATBF}). It has been adapted
to other situations.
For the case $R = \Z$, its role is taken by \cite[Theorem~10]{Li:ECOS},
for Euclidean domains by \cite[Theorem~4.1]{KK:CAGB},
and in our presentation by Theorem~\ref{thm:Spols}. 
From the vast literature on Gr\"obner bases, we highlight the monographs
\cite{CL:IVAA, BW:GB, AL:AITG, GP:ASIT} and the survey paper
\cite{BK:GB}. The mathematical content of the present note
builds upon \cite{KK:CAGB,Li:ECOS}; our definition of $S$-polynomials
differs from the one given in~\cite{Li:ECOS}, was inspired
by \cite{Bu:ACPA} and is close to \cite[Definition CP1]{KK:CAGB}.
The proof of Theorem~\ref{thm:Spols} was
modelled after the proof of Theorem 2.3.10 in \cite{Sm:ITAG},
and the notation using ``expressions with remainders'' follows~\cite{Ei:CA}.
I am indebted
to M.\ Kauers for sharing the unpublished notes
for his course on Gr\"obner bases at JKU in 2011 with me.
The presentation of how to solve linear systems in Theorem~\ref{sat:linsys}
was inspired by his lectures on linear algebra. 

The goal of the present note is to provide
a self-contained presentation
of as much of Gr\"obner basis theory as is needed to solve
linear systems of over $\Z[x_1, \ldots, x_n]$ (or $\rxn$ for a Euclidean domain $R$) in 
Section~\ref{sec:la}, along with proofs that are ready for
the classroom.
Some facts on partial orders go beyond
the material one commonly presupposes in an undergraduate
course. These facts are collected
in Section~\ref{sec:order}.
Most of the theory contained in the present note is
well known, but has so far been scattered in various research publications,
sometimes also with variations in the definitions. We
aim at providing one coherent presentation of these beautiful methods.

\section{Basic definitions}
We write $\N$ for the set of positive integers and $\ul{k}$ for
the set $\{1,\ldots, k\}$. For a set $G$, we write  ${G \choose 2}$
for
the set of two-element subsets of $G$.
\begin{de}
  Let $R$ be an integral domain. $R$ is a
  \emph{Euclidean domain} if there is a well ordered set
  $W$ and a map
  $\delta : R \to W$ such that
  $\delta(0) \le \delta(r)$ for all $r \in R$,
  and  for all $a, b \in R$ with $a \neq 0$, we have
  that
  $\delta(b) \le \delta(ab)$ and that there exist $q, r \in R$
  such that $b = qa + r$ and $\delta (r) < \delta (a)$.
\end{de}
Let $R$ be a Euclidean domain. For the polynomial ring
$\rvx := \rxn$, we will consider its module $\rvxk$, and
we will call the elements of this module \emph{polynomial vectors}.
For $p \in \rxn$ and $i \in \ul{k}$,
we write $p \, e_i$ for the vector $(0,\ldots, 0, p, 0,\ldots, 0)$
with $p$ at place $i$.
For $\alpha \in \NON$, we write $\mono{\alpha}$
for $x_1^{\alpha_1} \cdots x_n^{\alpha_n}$.
The elements of
$\{ \term{a}{\alpha}{i}  \mid a \in R \setminus\{0\},
                              (\alpha, i) \in \NONk \}$
are called \emph{term vectors} and the elements of
$$\monosnk := \{\term{}{\alpha}{i}  \mid 
(\alpha, i) \in \NONk \}$$
are called \emph{monomial vectors}.
We say that 
the monomial vector $\term{}{\alpha}{i}$ \emph{divides}
the monomial vector $\term{}{\beta}{j}$ and write
$\monov{\alpha}{i} \mid \monov{\beta}{j}$ if
$i = j$ and $\alpha_m \le \beta_m$ for all $m \in \ul{n}$.
This holds if and only if there is a monomial $\mono{\gamma}$ such that
$\mono{\gamma} \term{}{\alpha}{i} = \term{}{\beta}{j}$.
In this case, we will also write
$\frac{\term{}{\alpha}{i}}{\term{}{\beta}{j}}$ for $\mono{\gamma}$.
We say that the term vector $s = \term{a}{\alpha}{i}$ \emph{divides} the
term vector $t = \term{b}{\beta}{j}$ if $a$ divides $b$ in $R$
and $\term{}{\alpha}{i}$ divides $\term{}{\beta}{j}$.
In this case we write
$s \mid t$, and we write $\frac{t}{s}$ for the term $q \mono{\alpha}$ with
$(q \mono{\alpha})s=t$.
We can write every element from $\rvxk$ as a sum
$\sum_{(\alpha, i) \in E} \term{c_{(\alpha, i)}}{\alpha}{i}$, where $E$
is a finite subset of $\NONk$. A fundamental
step working with these polynomial vectors is to \emph{order} the terms in this sum
so that we can speak about a \emph{leading term vector}. A requirement
for this ordering is that when $\sum_{j=1}^m t_j$ is a sum of terms
that are in decreasing order, then so is the sum
$\sum_{j=1}^m \mono{\alpha}t_j$ obtained by multiplying with
a monomial $\mono{\alpha}$. This is reflected in Condition~\eqref{it:add}
of the following definition. Condition~\eqref{it:div} expresses
the requirement
that for all monomial vectors $m_1, m_2$ with $m_1 \mid m_2$,
we have $m_1 \le m_2$.
  \begin{de}
  Let $n, k \in \N$, and let $\lea$ be an order on $\monosnk$. This
  order $\le$ is \emph{admissible} if
  \begin{enumerate}
  \item $\lea$ is a total ordering;
  \item \label{it:div} for all $\term{}{\alpha}{i}, \term{}{\beta}{j} \in
    \monosnk$ with
    $\term{}{\alpha}{i} \mid \term{}{\beta}{j}$, we have
    $\term{}{\alpha}{i} \lea \term{}{\beta}{j}$;
\item \label{it:add}  for all $\alpha, \beta, \gamma \in \NON$
    and for all $i,j \in \ul{k}$ with
    $\term{}{\alpha}{i} \lea \term{}{\beta}{j}$, we have
    $\term{}{\alpha + \gamma}{i} \lea \term{}{\beta + \gamma}{j}$.
  \end{enumerate}   
  \end{de}
  An ordering $\lea'$ on $\NONk$ is admissible if
  the ordering $\lea$ defined by
  $\monov{\alpha}{i} \lea \monov{\beta}{j} :\Leftrightarrow
  (\alpha, i) \lea' (\beta, j)$ is admissible.

  One such ordering is the \emph{lexicographic position over term ordering},
  where for two distinct $(\alpha, i)$ and $(\beta, j) \in \NONk$, we have
  $((\alpha_1, \ldots, \alpha_n), i) <_{\textrm{lex}} ((\beta_1, \ldots, \beta_n), j)$ if
  $i > j$ or ($i = j$ and $\alpha_l < \beta_l$ for
  $l := \min \{m \in \ul{k} \mid \alpha_m \neq \beta_m\}$).
  Other admissible orderings can be defined by
  choosing a matrix $U \in \R^{n' \times n}$ and a permutation
  $\pi$ of $\ul{k}$ such that
  $\{\gamma \in \Q^n \mid U\gamma = 0\} = \{0\}$ and
  the first nonzero entry in every column of $U$ is positive.
  Then one can define an admissible order $\le_{U,\pi}$ by
  $(\alpha, i) \le_{U,\pi} (\beta, j) : \Leftrightarrow
   (U\alpha, \pi(i)) \le_{\textrm{lex}} (U\beta, \pi(j))$.

 For a finite subset $E$ of $\NONk$, an admissible ordering
$\lea$ of $\NONk$, a function
$c : E \to R$, and an $f \in \rxnk$ with $f \neq 0$ given by 
\(
f = \sum_{(\alpha, i) \in E} \term{c_{(\alpha, i)}}{\alpha}{i},
\)  
we define
\[
   \md (f) := {\max}_{\le} \{(\alpha,i) \in \NONk \mid c_{(\alpha, i)}
\neq 0 \};
\]
$\md (0)$ is not defined. 

Suppose that $f \neq 0$ and $(\gamma, i) = \md(f)$. Then we define
\[
\LM (f) := \term{}{\gamma}{i},  \,\, \LC (f) := c_{(\gamma, i)},
\,\, \LT (f) := \term{c_{(\gamma, i)}}{\gamma}{i}
\]
and call them the \emph{leading monomial vector}, the \emph{leading coefficient} and
the \emph{leading term vector}, respectively. All of these are undefined for $f=0$.
An important fact is that an admissible ordering on $\NONk$ is a
\emph{well order}, i.e., it is total and has no infinite strictly descending chains.
A proof is given in Lemma~\ref{lem:adwell}. This also implies
that for every nonempty subset $I$ of $\rvxk \setminus \{0\}$, there is
at least one $f \in I$ such that there is no $g \in I$
with $\md (g) < \md (f)$.
\begin{de} \label{de:sgb}
  Let $R$ be a Euclidean domain, let $I$ be a submodule
  of $\rxnk$, and let $\le$ be an admissible order of
  the monomial vectors.
  Then $G \subseteq I \setminus \{0\}$ is a \emph{strong
  Gr\"obner basis} of $I$ with respect to $\le$ if and only if
  for every $f \in I \setminus \{0\}$, there is an element $g \in G$
  such that $\LT(g) \mid \LT (f)$.
\end{de}
We write $\submod{G}$ for the submodule of $\rvxk$ generated by $G$.
When $G$ is a strong Gr\"obner basis of $I$, then
$\submod{G} = I$: Suppose that $\submod{G}$ is a proper subset of $I$,
and let $f$
be a polynomial vector of minimal degree $\md (f)$ in $I \setminus \submod{G}$
with respect to the admissible
ordering $\le$. The existence of such an $f$ -- under the assumption
$\submod{G} \neq I$ -- is justified by
Lemma~\ref{lem:adwell}.
Taking $g \in G$ such that $\LT (g) \mid \LT (f)$, we compute
$f' := f - \frac{\LT(f)}{\LT(g)} g$. If $f'=0$, then
$f = \frac{\LT(f)}{\LT(g)} g$ lies in $\submod{G}$.
If $f' \neq 0$, then by minimality, $f'$ lies in
$\submod{G}$, and hence so does $f = f' + \frac{\LT(f)}{\LT(g)} g$,
a contradiction.

When $R$ is a Euclidean domain, $\delta : R \to W$ is the grading
function of $R$ and
$f \in \rxnk \setminus \{0\}$,
we define the \emph{degree with $\delta$} of $f$ by
\[
\DegD (f) := (\md (f), \, \delta (\LC (f)));
\]
hence $\DegD (f) \in (\NONk) \times W$.
For a subset $I$ of $\rvxk$, we define
$\DegD (I) := \{ \DegD (f) \mid f \in I \setminus \{0\} \}$.
On $\NONk \times W$ we define an order
by
  \begin{equation} \label{eq:lmd}
 (\monoexp{\alpha}{i}, d) \divd (\monoexp{\beta}{j}, e) :\Longleftrightarrow
   \mono{\alpha} \mid \mono{\beta}, \, i = j, \, d \le e.
  \end{equation}
    Then the ordered set
  $(\NONk \times W, \divd)$ is isomorphic to the direct product of
  $n$ copies of $(\N_0, \le)$ with
  $(\{1,\ldots, k\}, =)$ and $(W, \le)$.
  Therefore the ordered set $(\NONk \times W, \divd)$
  has no infinite descending chains and 
  no infinite antichains (Theorem~\ref{thm:divdorder}\eqref{it:d1}).

\section{Existence of strong Gr\"obner bases}
\begin{thm} \label{thm:existssgb}
  Let $R$ be a Euclidean domain, let $I$ be a submodule
  of $\rxnk$,
  and let $\le$ be an admissible order of
  the monomial vectors.
  Then $I$ has a finite strong Gr\"obner basis with respect to $\le$.
\end{thm}
\begin{proof}
  Let $\Min (\DegD (I))$ be the set of minimal elements
  of $\DegD(I)$ with respect to the ordering $\divd$. Since
  $\Min (\DegD (I))$ is an antichain of $(\NONk \times W, \divd)$,
  it is finite (cf. Theorem~\ref{thm:divdorder}\eqref{it:d2}).
  Let $G$ be a finite subset of $I$ such that for every
  $(\monoexp{\alpha}{i}, d) \in \Min (\DegD (I))$, there
  is a $g \in G$ with
  $(\md (g),\delta (\LC (g))) = (\monoexp{\alpha}{i}, d)$.

  We claim that $G$ is a strong Gr\"obner basis.
  To show this, let $f \in I \setminus \{0\}$. Since
  $(\md (f), \delta (\LC(f))) \in \DegD (I)$, there is
  an $(\monoexp{\alpha}{i}, d) \in \Min(\DegD (I))$ with
  $(\monoexp{\alpha}{i}, d) \divd (\md (f), \delta (\LC(f)))$.
  Hence
  \[
  L := \{ g \in G \mid \DegD (g) \divd \DegD (f) \}
  \]
  is not empty. Let $g_1$ be an element of $L$ for which
  $\delta (\LC (g_1))$ is minimal. Since $\DegD(g_1) \divd \DegD(f)$,
  $\LM(g_1)$ divides $\LM (f)$. By the Euclidean property,
  there are $q, r \in R$ such that
  $\LC (f) = q \,\LC(g_1) + r$ with $\delta (r) < \delta (\LC(g_1))$.
  If $r = 0$, then $\LC(g_1) \mid \LC (f)$ and therefore
  $\LT (g_1)$ divides $\LT (f)$. Then $g_1$ is the required element
  from $G$.
  If $r \neq 0$, we let
  \[
  h := f - q \frac{\LM(f)}{\LM(g_1)} g_1.
  \]
  Then $h \in I$ and $\LT (h) = r \, \LM (f)$.
  Then there is $g_2 \in G$ such that
  $\DegD (g_2) \divd \DegD (h)$. Hence $\LM (g_2) \mid \LM (h)$ and
  $\delta (\LC (g_2)) \le \delta (r)$, and thus
  $\delta (\LC (g_2)) < \delta (\LC (g_1))$.
  Since $\delta (\LC (g_2)) < \delta (\LC (g_1)) \le \delta (\LC (f))$,
  we have $g_2 \in L$. This $g_2$ contradicts the minimality
  of $\LC(g_1)$. Therefore the case $r \neq 0$ cannot occur.
\end{proof}

\section{A criterion for being a strong Gr\"obner basis} \label{sec:crit}
In this section, we prove a criterion (Theorem~\ref{thm:Spols})
that guarantees that certain sets are strong Gr\"obner bases.
This criterion is then fundamental for constructing these bases
in Section~\ref{sec:construct}.
Throughout Sections~\ref{sec:crit} and~\ref{sec:construct}, $R$
will denote a Euclidean domain with grading function $\delta$.
We first need a generalization of Euclidean division, i.e., of
expressing $b$ as $qa +r$ with $\delta(r) < \delta (a)$, from
$R$ to $\rvx$.
\begin{de}
Let $G \subseteq \rvxk \setminus \{0\}$, and let $f \in \rvxk$. We say that
$\rho = ((a_i, m_i, g_i)_{i \in \ul{N}}, r)$ is an \emph{expression}
of $f$ by $G$ with remainder $r$ 
if $N \in \N_0$ and  for each $i \in \ul{N}$, we have that
$a_i \in R$, $m_i$ is a monomial,
$g_i \in G$, $r \in \rvxk$ and
\[
f = \sum_{i=1}^N a_i m_i g_i + r.
\]
An expression is \emph{Euclidean} with respect to
the admissible monomial vector ordering $\lea$ if
 for all $i \in \ul{N}$, we have 
  $\LM (m_i g_i) \lea \LM (f)$,
             and ($r=0$ or there is no
             $g \in G$ such that
             $\DegD (g) \divd \DegD (r)$).
\end{de}
We note that in an expression, $a_i = 0$ is allowed.
 The name \emph{expression} follows the notation of
 \cite[Definition~15.6]{Ei:CA}.
We will construct such expressions using Euclidean division.
\begin{alg}[Euclidean division] \label{alg:ed1}
       \mbox{} \\
       Input: $f \in \rvxk \setminus \{0\}$, $G \subseteq \rvxk
       \setminus \{0\}$, an admissible order $\le$ of $\NONk$.
       \\
       Ouput: An Euclidean expression $((a_i, m_i, g_i)_{i \in \ul{N}}, r)$
       of $f$ by $G$.
             \begin{algorithmic}[1]
               \State $r \gets f$
               \State $\rho \gets ()$
               \While{$r \neq 0$ and $\exists g \in G : \DegD(g) \divd \DegD (r)$}
               \State Find some $q, s \in R$ with $\LC(r) = q \, \LC(g) + s$
               and $\delta (s) < \delta (\LC(g))$.
               \State $r \gets r - q \frac{\LM(r)}{\LM(g)} g$
               \State Append $(q, \frac{\LM(r)}{\LM(g)}, g)$ to $\rho$
               \EndWhile
             \State Return $(\rho, r)$  
         \end{algorithmic}
\end{alg}
\begin{lem}
   For each input $f, G$,
   Algorithm~\ref{alg:ed1} terminates and yields a Euclidean expression of $f$ by $G$.
\end{lem}
\begin{proof}
  We first prove termination.
  We say that $f \in \rvxk$ does not guarantee termination
  if there is an infinite sequence  $r_0 = f, r_1, r_2\ldots$
  of values of $r$ produced by the algorithm.
  Among those $f$ that do not guarantee termination, we let
  $L$ be the set of those $f$ for which $\LM (f)$ is minimal with respect
  to the admissible ordering $\lea$. Among the elements of $L$, we choose
  $f$ to be of minimal $\delta (\LC(f))$.
  If in the computation $r_1 = f - q \,\frac{\LM(f)}{\LM(g)} g$,
  with $s = \LC(f) - q \, \LC(g)$, we have $s = 0$,
  then $\LM (r_1) \strictlea \LM(f)$. Then $r_1$ does not
  guarantee termination, contradicting the minimality of
  $\LM(f)$.
  If $s \neq 0$ and $\delta (s) < \delta (\LC(g))$,
  we have $\LT (r_1) = s \,\LM (f)$ and therefore
  $\delta (\LC(r_1)) = \delta (s) < \delta (\LC(g))$.
  Since $\DegD (g) \divd \DegD (f)$, we have
  $\delta (\LC(g)) \le \delta (\LC(f))$.
  Thus $\delta (\LC(r_1)) < \delta (\LC (f))$.
  Since $r_1$ does not guarantee termination,
  we have a contadiction to
  the minimality of $\delta (\LC (f))$.

  For proving correctness, we observe that
  throughout the algorithm $(\rho, r)$ is an expression of $f$ by $G$
  satisfying the degree bound. When the while-loop is left,
  then $r$ has the required properties.
\end{proof}  

 Expressions with remainder~$0$
 will also be called \emph{representations}. The importance of
 representations in which only one summand has maximal degree
 was observed in \cite{Li:ECOS}.
 \begin{de}
   Let $f \in \rvxk \setminus \{0\}$ and $G \subseteq \rvxk\setminus \{0\}$.
   Then $\rho = (a_i, m_i, g_i)_{i \in \ul{N}}$ is a 
  \emph{strong standard representation} of $f$ by $G$ with respect to
the monomial vector ordering $\lea$ if $(\rho, 0)$ is an expression
of $f$ by $G$ with remainder~$0$,
and in addition,
\[
 N \ge 1, \, \LM (m_1 g_1) = \LM (f), \text{ and }
\LM (m_i g_i) \strictlea \LM(f) \text{ for all } i \in \ul{N} \setminus \{1\}.
\]
 \end{de}
 
 We will now define $S$-polynomial vectors; this definition is
 a slight modification of \cite[Definition CP1]{KK:CAGB}.
For $\alpha, \beta \in \NON$, we let
$\alpha \sqcup \beta := (\max (\alpha_1, \beta_1), \ldots, \max (\alpha_n, \beta_n))$.
Hence $\mono{\alpha\sqcup\beta}$ is the least common multiple
of $\mono{\alpha}$ and $\mono{\beta}$ in $\rvx$.
\begin{de}[$S$-polynomial vectors] \label{de:spol}
Let $f,g \in \rvxk \setminus \{0\}$ with $f \neq g$, and 
assume that $\LT (f) = \term{a}{\alpha}{i}$ and
$\LT (g) = \term{b}{\beta}{j}$. Let
\[ \alpha' :=
(\alpha \sqcup \beta) - \beta \text{ and } \beta' := (\alpha \sqcup \beta) - \alpha.
\]
Then $h \in \rvxk$ is an {\emph{$S$-polynomial}} \emph{vector}
of the pair  $(f,g)$ if one of the following two conditions holds:
\begin{enumerate}
    \item \label{it:SC1}
$i = j$, $\delta(a) \ge \delta(b)$ and there exists $q \in R$
such that $\delta(a - q b) < \delta(a)$ and
\[
h = \mono{\beta'} f - q \mono{\alpha'} g;
\]
  \item
    $i \neq j$ and $h = 0$.
\end{enumerate}    
The polynomial vector $h$ is an $S$-polynomial vector
of the set $\{f,g\}$ if $h$ is an $S$-polynomial vector
of $(f,g)$ or of $(g,f)$.
\end{de}
Concerning item~\eqref{it:SC1}, we notice that Euclidean
division of $a$ by $b$ in $R$ would yield a $q$ that even satisfies
$\delta (a - qb) < \delta (b)$, but for our purposes the weaker
condition
$\delta (a - qb) < \delta (a)$ suffices.
\begin{thm} \label{thm:Spols}
  Let $\le$ be an admissible ordering on $\monosnk$,
  and let $G \subseteq \rvxk \setminus \{0\}$. We assume that for all $f,g \in G$
  with $f \neq g$,
  there is an $S$-polynomial vector $h$ of $\{f,g\}$ such that $h = 0$
  or $h$
  has a strong standard representation. %
  Then $G$ is a strong Gr\"obner basis with respect to $\le$
  for the submodule $I$ of $\rvxk$ that is
  generated by $G$. 
\end{thm}
\begin{proof}
  Let $f \in I \setminus \{0\}$. We will show that there
  is $g \in G$ with $\LT(g) \mid \LT(f)$.
  Since $f \in I$, there is
  $\rho = (a_i, m_i, g_i)_{i \in \ul{N}}$ such that
  $f = \sum_{i \in \ul{N}} a_i m_i g_i$.
  Such a $\rho$ is called a \emph{representation} of $f$ by $G$.
  Here, no restriction on $\md (m_i g_i)$ is made.
   We measure the complexity of a
   representation $\rho = (a_i, m_i, g_i)_{i \in \ul{N}}$
   using the following complexity parameters:
   \begin{equation*}
           C_1 (\rho) :=  \max \, \{\md (m_i g_i) \mid i \in
           \ul{N} \}
   \end{equation*}
   is the maximal degree of $m_i g_i$ appearing in $\rho$,
   where the maximum is taken with respect to the admissible
   ordering of on $\NONk$. We let
   \[
   I_1 (\rho)  =  \{ i \in \ul{N} \mid \md (m_i g_i) = C_1 (\rho) \}
   \]
   be the set of those indices for which this maximum is attained.
   We define
   \[ 
       C_2 (\rho) :=  \max \, \{\delta (\LC(g_i)) \mid
       i \in I_1 (\rho) \}
   \]
   as the maximum of the $\delta$-grades of the leading coefficients
   of those $g_i$'s for which $\md (m_i g_i)$ is maximal.
   The set 
   \[
        I_2 (\rho) :=  \{ i \in \ul{N} \mid
         \md (m_i g_i) = C_1 (\rho) \text{ and }
         \delta (\LC(g_i)) = C_2 (\rho) \}
   \]
   collects those indices from $I_1 (\rho)$ for which this maximum of $\delta$-grades
   are attained. Finally,
   \[
   C_3 (\rho) := \# I_2
   \]
   counts the number of elements of $I_2$.
   Now we choose a representation $\rho = (a_i, m_i, g_i)_{i \in \ul{N}}$
   of $f$ for 
  which the triple
  $(C_1 (\rho),C_2 (\rho), C_3 (\rho))$
  is minimal with respect
  to the lexicographic ordering on $(\NONk) \times W \times \N$, where
  the order on $\NONk$ is taken to be the admissible ordering $\le$.
  This means that $\rho$ minimizes
   $C_1$ with respect to the admissible order on $\NONk$,
   among those that minimize~$C_1$, $\rho$
   minimizes $C_2$, and so on.
   Since all three sets $\NONk,
   W, \N$ are well ordered, i.e., totally ordered without
   infinite descending chains, such
   a minimizing $\rho$ exists.
  Since for every permutation of $\ul{N}$, the
  representation $\rho' = (a_{\pi(i)}, m_{\pi(i)}, g_{\pi(i)})_{i \in \ul{N}}$ of $f$ has
  the same complexity parameters
  as $\rho$, we may assume
   \begin{equation} \label{eq:degreeo}
  \md (m_1 g_1) \gea \md  (m_2 g_2) \gea  \cdots \gea
  \md (m_N g_N).
   \end{equation}
  We now consider several cases:

  \textbf{Case 1}: $\#I_1 (\rho) = 1$:
   By the assumption~\eqref{eq:degreeo}, we then have $I_1 (\rho) = \{1\}$.
    If $a_1 = 0$, we take the
  representation $\rho' := (a_i, m_i, g_i)_{i \in \ul{N} \setminus \{1\}}$. Then
    $C_1 (\rho') < C_1 (\rho)$, contradicting the
  minimality of $\rho$.
  If $a_1 \neq 0$, then
  $\md (a_1 m_1 g_1) = \md (f)$ and
  $\LC (f) = a_1 \, \LC (g_1)$. Therefore
  $\LT (g_1) \mid \LT (f)$. 

    \textbf{Case 2}: $\#I_1 (\rho) \ge 2$:
  Let $l := \# I_1 (\rho)$.
  Then $\LM (m_1 g_1) = \cdots = \LM (m_{l} g_{l})$, and we
  may assume without loss of generality that 
  \begin{equation} \label{eq:deltao}
  \delta (\LC(g_1)) \ge \delta (\LC(g_2)) \ge \cdots \ge \delta(\LC(g_{l})).
  \end{equation}

  \textbf{Case 2.1}: $g_1 = g_2$: Since $\LM (m_1 g_1) = \LM (m_2 g_2)$,
      we then have $m_1 = m_2$. Hence
  $\rho' := ((a_1 + a_2, m_2, g_2), (a_3, m_3, g_3), \ldots,
  (a_N, m_N, g_N))$  is a representation of $f$
      with $C_1 (\rho') = C_1 (\rho)$. Since
      $\delta (\LC(g_1)) = \delta (\LC(g_2))$, we also have
   $C_2 (\rho') = C_2 (\rho)$. Now 
  $C_3 (\rho') = C_3 (\rho) - 1 < C_3 (\rho)$.
  Then $\rho'$ contradicts the
  minimality of $\rho$.

  \textbf{Case 2.2}: $g_1 \neq g_2$:
   Let
   \[
   \LT (g_1) = \term{a}{\alpha}{i} \text{ and }
   \LT(g_2) = \term{b}{\beta}{j}.
   \]
  Since $\md (m_1 g_1) = \md (m_2 g_2)$, we have $i = j$.
  By the assumptions, the $S$-polynomial vector $h$ coming from
  $\{g_1, g_2\}$ is $0$ or has a strong standard representation by $G$.
    If $\delta(g_1) > \delta(g_2)$, then $h$ is an $S$-polynomial vector
  of the pair $(g_1, g_2)$.
  If $\delta (\LC(g_1)) = \delta (\LC(g_2))$
  and $h$ is an $S$-polynomial vector of the pair $(g_2, g_1)$,
  we swap the first two entries in the representation $\rho$ and
  obtain a representation $\tilde{\rho}$ that still satisfies~\eqref{eq:degreeo}
  and~\eqref{eq:deltao}. This allows us to assume that $h$ is
  an $S$-polynomial vector of the pair $(g_1, g_2)$.
   Let
   \[ \alpha' := (\alpha \sqcup \beta) - \beta, \,\,
   \beta' := (\alpha \sqcup \beta) - \alpha,
   \]
   and let 
  $\gamma$ be such that $\gamma + (\alpha \sqcup \beta) =
  \md (m_1 g_1)$.
  Let $q \in R$ be such that $\delta (a-qb) < \delta (a)$
  and 
  \[
  h = \mono{\beta'} g_1 - q \mono{\alpha'} g_2.
  \]
  Then
  \[
  \mono{\gamma} h = m_1 g_1 - q m_2 g_2,
  \]
  and thus
  \begin{equation} \label{eq:star}
  m_1 g_1 = \mono{\gamma} h + q m_2 g_2.
  \end{equation}

  \textbf{Case 2.2.1}: $h = 0$: In this case,
  $m_1 g_1 = q m_2 g_2$, and thus
   \[
   \rho' := ((a_1 q + a_2, m_2, g_2), (a_3, m_3, g_3), \ldots, (a_N, m_N, g_N))
   \]
   is a representation of $f$ that
  satisfies $C_1 (\rho') = C_1 (\rho)$.

  \textbf{Case 2.2.1.1}: $\delta (\LC(g_1)) > \delta (\LC (g_2))$:
  Then 
  $C_2 (\rho') = \delta ( \LC(g_2) ) < \delta (\LC(g_1)) = C_1 (\rho')$.
  Thus $\rho'$ contradicts the minimality of $\rho$.

 \textbf{Case 2.2.1.2}: $\delta (\LC(g_1)) = \delta (\LC (g_2))$:
 Then $C_2 (\rho') = \delta (\LC(g_2))$ and
 $C_3 (\rho') = C_3 (\rho) - 1$, contradicting the minimality of $\rho$.
 
 \textbf{Case 2.2.2}: $h \neq 0$:
  By the assumptions,
  $h$ has a strong standard representation 
  $(b_i, n_i, h_i)_{i \in \ul{M}}$ with the $h_i$'s in $G$.
  Now from~\eqref{eq:star}, we obtain
  \[
  a_1 m_1 g_1 = \sum_{i \in \ul{M}}
  a_1 b_i (\mono{\gamma} n_i) h_i + a_1 q m_2 g_2,
  \]
  and therefore
  \[
  a_1 m_1 g_1 + a_2 m_2 g_2 =
  \sum_{i \in \ul{M}}
  a_1 b_i (\mono{\gamma} n_i) h_i + (a_1 q + a_2) m_2 g_2.
  \]
  We claim that the representation $\rho'$ coming from
  \[
     f =  \big( \sum_{i \in \ul{M}}
     a_1 b_i (\mono{\gamma} n_i)  h_i \big) + (a_1 q + a_2) m_2 g_2 +
     \sum_{i=3}^N a_i m_i g_i
     \]
     has lower complexity than $\rho$.
     We know that 
     $\md (h) \le \md (\mono{\beta'} g_1)$ and
     thus $\md (\mono{\gamma} h) \le \md (m_1 g_1)$. We distinguish
     cases according to whether this inequality is strict.
     
  \textbf{Case 2.2.2.1}: $\md (\mono{\gamma} h) < \md (m_1 g_1)$:
  Then for all $i \in \ul{M}$, we have
  $\md (\mono{\gamma} n_i h_i) < \md (m_1 g_1)$.
  Since $\md(m_1 g_1) = \md (m_2 g_2)$, we therefore have
  $C_1 (\rho') = C_1 (\rho)$.
  If 
  $\delta (\LC (g_1)) > \delta (\LC (g_2))$, we have 
  $C_2 (\rho') < C_2 (\rho)$, and
  if $\delta (\LC (g_1)) = \delta (\LC (g_2))$, we have
  $C_2 (\rho') = C_2 (\rho)$ and $C_3 (\rho') = C_3 (\rho) - 1 <
   C_3 (\rho)$, contradicting the minimality of $\rho$.

  \textbf{Case 2.2.2.2}: $\md (\mono{\gamma} h) = \md (m_1 g_1)$:
  Then $\md (\mono{\gamma} n_1  h_1) =\md (m_1 g_1)$
  and
  $\md (\mono{\gamma} n_i  h_i) < \md (m_1 g_1)$
  for $i \in \{2, \ldots, M\}$.
  Hence $C_1 (\rho') = C_1 (\rho)$.
  Since $\LT (b_1 n_1 h_1) = \LT (h)$,
  we have $\delta (\LC(h_1)) \le \delta (b_1 \, \LC (h_1))
  = \delta (\LC (h))$.
  From the definition of $S$-polynomial vectors, we
  have $\delta (\LC(h)) < \delta (\LC (g_1))$, and thus
  $\delta (\LC(h_1)) < \delta (\LC(g_1))$.
  If 
  $\delta (\LC (g_1)) > \delta (\LC (g_2))$, we have 
  $C_2 (\rho') < C_2 (\rho)$, and
  if $\delta (\LC (g_1)) = \delta (\LC (g_2))$, we have
  $C_2 (\rho') = C_2 (\rho)$ and $C_3 (\rho') = C_3 (\rho) - 1 <
  C_3 (\rho)$, contradicting the minimality of $\rho$.

Hence in Case~2, we always obtain $\rho'$ with complexity
           than $\rho$, showing that the
           case $\#I_1 (\rho) \ge 2$ cannot occur.
\end{proof}

\section{Construction of strong Gr\"obner bases} \label{sec:construct}
In this section, we assume that a submodule $I$ of $\rvxk$ is
given by a finite set $F$ of generators. (By Hilbert's Basis Theorem, or
simply by Theorem~\ref{thm:existssgb}, such a finite $F$ exists.)
Our goal is to construct a finite strong Gr\"obner basis $G$ for $I=\submod{F}$.
We will proceed by  adding polynomials to $F$
in order to obtain a set $G$ such that each
$2$-element subset of $G$ has an $S$-polynomial vector
with a strong standard representation; then Theorem~\ref{thm:Spols} guarantees
that we have found a strong Gr\"obner basis.
In one step, we consider one $2$-element subset $\{p,q\}$ of~$F$.
The augmentation of $F$ using the set $\{p,q\}$ yields a set $F'$
in which either $\{p,q\}$ has an $S$-polynomial vector that has a strong
standard representation, or $\{p,q\}$ still has no strong standard representation
in $F'$, but $F'$ is, in some sense, \emph{larger} than $F$, which also brings us closer
to termination.
To express this idea of becoming larger, we let $W$ be the codomain of the
  Euclidean grading function $\delta$. We consider subsets $T$
  of $\NONk \times W$ and we say that such a subset $T$ is
  \emph{upward closed} if for all $s \in T$ and $t \in \NONk \times W$
  with $s \divd t$, 
  we have $t \in T$. Since $(\NONk \times W, \divd)$ has no
  infinite descending chains and no infinite antichains, there
  is no infinite ascending chain of upward closed subsets
  of $\NONk \times W$ with respect to $\subseteq$ (cf. Theorem~\ref{thm:divdorder}\eqref{it:d3}).
  For $G \subseteq \rvxk$, we will consider the upward closed set
  \[
  \ups{\DegD (G)} := \{ (\monoexp{\gamma}{i}, d) \in \NONk \times W \mid
  \exists g \in G : \DegD(g) \divd (\monoexp{\gamma}{i}, d) \}.
  \]

  \begin{alg}[Augmentation] \label{alg:augment} \mbox{} \\
  Input: A finite subset $G$ of $\rvxk$, a two element
  subset $\{p, q\}$ of $G$, and an admissible order $\le$ of $\NONk$. \\
  Output: A pair $(H, x)$, where
  $H$ is a finite set with $G \subseteq H \subseteq \submod{G}$
    and $x \in \{0,1\}$ such that the following hold:
  \begin{enumerate}
  \item \label{it:aug1}  If $x = 1$, then $\{p, q\}$ has an $S$-polynomial vector $f$ such  that
    $f = 0$ or $f$  has a strong standard representation by $H$.
  \item  \label{it:aug2} If $G \neq H$, then $\ups{\DegD (G)} \subset \ups{\DegD (H)}$.
  \item \label{it:aug3} If $G = H$, then $x = 1$.   
  \end{enumerate}  
  \begin{algorithmic}[1]
    \Function{Augment}{$G, \{p, q\}$}
    \State $f \gets$ some $S$-polynomial vector of $\{p, q\}$ \label{li:e1}
    \State $x \gets 0$
         \If{$f = 0$}
             \State $x \gets 1$
         \ElsIf{$\exists g \in G : \LT (g) \mid \LT (f)$}
            \State $x \gets 1$ \label{line:P}
            \State $f' \gets f - \frac{\LT(f)}{\LT(g)} g$
            \State Find a Euclidean expression \label{li:e2}
             $f' = \sum_{i=1}^M a_i m_i g_i + r$ by $G$.
             \If{$r \neq 0$}
                 \State $G \gets G \cup \{r\}$ \label{li:assg1}
             \EndIf
        \ElsIf{$\exists g \in G : \DegD (g) \divd \DegD(f)$}
          \State 
            Among those $g \in G$ with
           $\DegD (g) \divd \DegD (f)$, pick $g$ with minimal
           $\delta (\LC (g))$.
           \State Find $q \in R$ such that
           $\delta \bigl(\LC(f) - q \, \LC(g)) < \delta (\LC(g)\bigr)$
           \State $f' \gets f - q \frac{\LM(f)}{\LM(g)} g$
           \State $G \gets G \cup \{f'\}$ \label{li:assg2}
        \Else
             \State $G \gets G \cup \{f\}$ \label{li:assg3}     
             \EndIf
         \State Return $(G,x)$    
    \EndFunction
   \end{algorithmic}     
 \end{alg}  

  \begin{lem} Algorithm~\ref{alg:augment} is correct.
  \end{lem}
  \begin{proof}
    For proving the first output condition, we assume that
    $x = 1$. We show that then $f$ is an $S$-polynomial vector of $\{p,q\}$
    with the required conditions. If $f=0$, this is clearly
    the case. If $\exists g \in G : \LT (g) \mid \LT(f)$,
    then  
  in the case $r = 0$,
  $((\frac{\LC(f)}{\LC(g)}, \frac{\LM(f)}{\LM(g)}, g),
  (a_1, m_1, g_1), \ldots, (a_Mm_Mg_M))$ is  a strong standard representation
  of $f$ by $G$
  and in the case
  $r \neq 0$, we note that $\LM(r) \lea \LM(f') < \LM(f)$ and thus
   $((\frac{\LC(f)}{\LC(g)}, \frac{\LM(f)}{\LM(g)}, g),
  (a_1, m_1, g_1), \ldots, (a_Mm_Mg_M), (1, \mono{0}, r))$
  is a strong standard representation. 

  For proving the second output condition, we assume $G \neq H$.
  If there exists $g \in G$ with $\LT(g) \mid \LT (f)$ and $r \neq 0$,
  then since $r$ is the remainder of a Euclidean division,
  $\DegD (r) \not\in \ups{\DegD (G)}$, and thus
  $\ups{\DegD (G)} \subset \ups{\DegD (G \cup \{r\})} = \ups{\DegD (H)}$.

  If there is no $g \in G$ with $\LT (g) \mid \LT (f)$, but 
  there is a $g \in G$ such that $\DegD (g) \divd \DegD (f)$,
    then we show that $\DegD (f') \not\in \ups{\DegD (G)}$.
    Suppose that there is $g_1 \in G$ with
    $\DegD (g_1) \divd \DegD (f')$. Then $\delta (\LC(g_1)) \le
     \delta (\LC(f'))$ and 
    $\LM (g_1) \mid \LM (f')$. Since $\LT (g) \nmid \LT (f)$,
    we have $\LM (f') = \LM(f)$, and thus $\LM(g_1) \mid \LM (f)$.
    Furthermore $\LC (f') = \LC (f) - q \, \LC (g)$ and thus
    $\delta (\LC(f')) < \delta (\LC(g))$, and from
    $\DegD (g) \divd \DegD (f)$ we obtain
    $\delta(\LC(g)) \le \delta (\LC(f))$.
    Altogether
    $\delta(\LC (g_1)) \le \delta(\LC (f')) < \delta(\LC (g)) \le
    \delta(\LC(f))$. From this, we obtain $\DegD(g_1) \divd \DegD(f)$
    and $\delta (\LC(g_1)) < \delta (\LC(g))$, contradicting the
    minimality of $\delta (\LC(g))$.
    Thus $\DegD (f') \not\in \ups{\DegD (G)}$, and 
    therefore
    $\ups{\DegD (G)} \subset \ups{\DegD (G \cup \{f'\})} = \ups{\DegD (H)}$.

    Finally, if $f \neq 0$ and
    $\exists g \in G : \DegD (g) \divd \DegD (f)$ is false, then
    $\DegD (f) \not\in \ups{\DegD (g)}$,
    and thus $\ups{\DegD (G)} \subset \ups{\DegD (G \cup \{f\})}
    = \ups{\DegD (H)}$.

    For proving the third output condition, we assume $G=H$.
    This happens only if $f=0$ or if
    ($f \neq 0$, $\exists g \in G : \LT (g) \mid \LT (f)$
    and $r = 0$) because in all other cases, a polynomial
    vector gets added to $G$. In both of these cases,
    $x$ is set to $1$.
 \end{proof}

\begin{alg}[Strong Gr\"obner Basis] \label{alg:sgb}
       \mbox{} \\
       Input: $F \subseteq \rvxk \setminus \{0\}$,
        an admissible order $\le$ of $\NONk$.
       \\
       Output: $G \subseteq \rvxk \setminus \{0\}$ such that
       is a strong Gr\"obner basis of $\submod{F}$
       with respect to the monomial vector ordering
       $\le$.
       \begin{algorithmic}[1]
         \State $G \gets F$.
         \State $P \gets \emptyset$.
         \While{$\exists p, q \in G : p \neq q $  and $\{p, q\} \not\in P$}
         \State $(G, x) \gets \textsc{Augment} (G, \{p,q\})$
         \If{$x=1$}
             \State $P \gets P \cup \{ \{p,q \} \}$
         \EndIf    
       \EndWhile
       \State Return $G$
   \end{algorithmic}
\end{alg}

\begin{thm} Algorithm~\ref{alg:sgb} terminates on every input and
  produces a correct result.
\end{thm}
\begin{proof}
  We first observe that throughout the algorithm,
  $G$ generates the same submodule as $F$.
  Furthermore, each $\{p,q\} \in P$ has an $S$-polynomial vector $f$
  that is $0$ or has a strong standard representation.
  The set $\{p,q\}$ can only be added to $P$ when $x=1$ and in this case the output
  condition~\eqref{it:aug1} of \textsc{Augment} guarantees that
  $\{p,q\}$ has $f$ as required.
Hence if the algorithm terminates, all two-element subsets
  of $G$ have a strong standard representation.
  Thus by Theorem~\ref{thm:Spols}, $G$ is then a strong Gr\"obner basis.

  In order to show termination, we let $W$ be the codomain of the
  Euclidean grading function $\delta$, and  we consider the upward closed
  subset 
    \[
  \ups{\DegD (G)} := \{ (\monoexp{\gamma}{i}, d) \in \NONk \times W \mid
  \exists g \in G : \DegD(g) \divd (\monoexp{\gamma}{i}, d) \}
  \]
  of $(\NONk \times W, \divd)$.
  Our claim is that in each execution of the while loop,
  if $G_1$ and $P_1$ are the values of $G$ and $P$ when entering
  the loop and $G_2$ and $P_2$ are the values before the next iteration
  of the loop, we have $\ups{\DegD (G_1)} \subset \ups{\DegD(G_2)}$
  or ($\ups{\DegD (G_1)} = \ups{\DegD(G_2)}$ and
  $\#(\VecTwo{G_2}{2} \setminus P_2) < \#(\VecTwo{G_1}{2} \setminus P_1)$\,). 
  If $G_1 \neq G_2$, then output condition~\eqref{it:aug2} of
  \textsc{Augment} yields $\ups{\DegD(G)} \subset \ups{\DegD(H)}$.
  If $G_1 = G_2$, then by output condition~\eqref{it:aug3} of
  \textsc{Augment}, we have $x=1$ and thus $P_2 = P_1 \cup \{p,q\}$
  and therefore
  $\#(\VecTwo{G_2}{2} \setminus P_2) < \#(\VecTwo{G_1}{2} \setminus P_1)$. 
    
    Now suppose that there is an execution of this algorithm that
        does not terminate. Then we know that from some
        point onwards, $\ups{\DegD(G)}$ stays constant, and from this
        point on, $\#(\VecTwo{G}{2} \setminus P)$
        strictly descends forever, which is impossible.
   \end{proof}

\cite[Theorem~11]{Li:ECOS} contains a
criterion\footnote{In the statement of \cite[Theorem~11]{Li:ECOS}, the
assumption $c_1 \in \{-1, +1\}$ is missing. Without adding this
assumption, for $p_1 := 2x + 1$ and $p_2 := 4 y + 1$, we obtain
$\operatorname{SPoly}_2 (p_1, p_2) = 2 y p_1 - x p_2 =
2 y- x$, which has no strong standard representation since
$2y$ and $x$ are not divisible by any of $2x$ and $4y$. --
In the proof given in \cite[Theorem~11]{Li:ECOS},
the $S$-polynomial of $p_1, q_1$ is computed (incorrectly)
as $\operatorname{SPoly}_2 (p_1, p_2) = 4 y p_1 - 2 x p_2 =
4y- 2x$.}
that  generalizes
\cite[p.377, S.2.]{Bu:EAKF}, which tells that certain $S$-polynomial
vectors need not be considered.
We provide a generalization, which, 
when dealing with polynomial
vectors, needs the rather restrictive assumption that
both polynomial vectors $f,g$ have entries only in the same component.
When speaking of polynomials instead of polynomial vectors,
we write the leading term of $p \neq 0$ as $\lt (p)$, the leading monomial as
$\lm (p)$, and the degree of $p$,
which is an element in $\N_0^n$, as
$\deg (p)$.
\begin{thm}
  Let $\tilde{f}, \tilde{g} \in \rvx, \alpha, \beta \in \NON$, $i \in \ul{k}$,
  and $a, u \in R$ such that $u$ is a unit in $R$,
  $\lt (\tilde{f}) = \termp{a}{\alpha}$ and
  $\lt (\tilde{g}) = \termp{u}{\beta}$. Let $f:= \tilde{f} \, e_i$ and
  $g := \tilde{g} \, e_i$. If $\termp{}{\alpha}$ and $\termp{}{\beta}$
  are coprime monomials (which means that for all $j \in \ul{n}$
  we have $\alpha_j = 0$ or $\beta_j = 0$),
    then $\{f ,g \}$ has an $S$-polynomial vector that is $0$
  or
  has a strong standard
  representation by $\{f,g\}$.
\end{thm}
\begin{proof}
  Let $\tilde{f}_1 := \tilde{f} - \lt (\tilde{f})$ and
  $\tilde{g}_1 := \tilde{g} - \lt (\tilde{g})$.
  Then $h = \mono{\beta} f  - au^{-1} \mono{\alpha} g$ is an
  $S$-polynomial vector of $\{f,g\}$. Suppose $h \neq 0$.
  We have
  \begin{multline} \label{eq:asrep}
  h = \mono{\beta} f - au^{-1} \mono{\alpha} g =
      (\mono{\beta} \tilde{f} - au^{-1} \mono{\alpha} \tilde{g}) \, e_i 
    =    
      u^{-1} (u\mono{\beta} \tilde{f} - a \mono{\alpha} \tilde{g}) \, e_i 
     \\
     =  u^{-1}
     \big( (\tilde{g} - \tilde{g}_1) \tilde{f} -
     (\tilde{f} - \tilde{f}_1) \tilde{g} \big) e_i
     =
     - u^{-1} \tilde{g}_1 f + u^{-1} \tilde{f}_1 g.
  \end{multline}
  Writing $\tilde{g}_1$ and $\tilde{f}_1$ as sums of terms,
  we obtain a representation of $h$.
  We will show now that it is a strong standard representation.
  If $f_1 = 0$ or $g_1 = 0$, then this representation has only one summand
  of degree $\md(h)$ and is therefore a strong standard representation.
  Hence let us assume $f_1 \neq 0$ and $g_1 \neq 0$.
  We observe that
  $\md (\lm (g_1) f) \neq \md (\lm (f_1) g)$:
  Seeking a contradiction, we assume
  $\md (\lm (g_1) f) = \md (\lm (f_1) g)$.
  Then $\lm (g_1) \lm(\tilde{f}) = \lm (f_1) \lm (\tilde{g})$,
  which means
  $\lm (g_1) \termp{}{\alpha} = \lm (f_1) \termp{}{\beta}$.
  We therefore have $\termp{}{\alpha} \mid \lm (f_1) \termp{}{\beta}$.
    By the assumptions on $\alpha$ and $\beta$, we then 
    have $\termp{}{\alpha} \mid \lm (f_1)$, contradicting
    $\deg (f_1) < \alpha$.

    Therefore, exactly one of
    $- u^{-1} \tilde{g}_1 f$ and $u^{-1} \tilde{f}_1 g$
    has degree $\md (- u^{-1} \tilde{g}_1 f +u^{-1} \tilde{f}_1 g)$,
    which is equal to $\md (h)$. Thus by writing $\tilde{g}_1$ and
    $\tilde{g}_1$ as sums of terms,~\eqref{eq:asrep} produces
    a strong standard representation
    of $h$ with respect to $\{f,g\}$.
\end{proof}

 \section{Existence and uniqueness of reduced strong Gr\"obner bases} \label{sec:redgb}
 The construction given in Section~\ref{sec:construct} has the shortcoming
 that during the process, polynomial vectors can never be removed
 from a basis. Also, once we have found a strong Gr\"obner basis $G$ of $I$,
 then we see from Definition~\ref{de:sgb} that every $G'$ with
 $G \subseteq G' \subseteq I \setminus \{0\}$ is also a strong Gr\"obner basis.
 Hence a strong Gr\"obner basis of $I$ need not be unique. However, we obtain
 uniqueness if we require that the Gr\"obner basis is \emph{reduced}. In this
 section, we prove the existence and uniqueness of such a reduced Gr\"obner basis;
 Section~\ref{sec:constructreduced} is then devoted to its
 algorithmic construction. 

 When $R$ is the Euclidean domain $\Z$ with grading function $\delta (z) := |z|$, then
 $6$ may be expressed by $4$ either as $6 = 1 \cdot 4 + 2$ or $6 = 2 \cdot 4 + (-2)$.
 In order to be able to prefer one of this expressions, we need to refine
 the grading function $\delta$. Hence when
 $R$ is a Euclidean domain with grading function $\delta : R \to W$, 
 we assume that we additionally have an injective function
 $\dhat$ from $R$ into a well ordered set $W'$ with the property
$\dhat (0) \le \dhat (a)$ for all $a \in R$ and 
\begin{equation} \label{eq:dhat}
 \text{ for all } a, b \in R : \delta (a) < \delta (b) \Rightarrow
                     \dhat (a) < \dhat (b).
\end{equation}
Throughout Sections~\ref{sec:redgb} and~\ref{sec:constructreduced},
we assume that $R$ is a Euclidean domain with the functions
$\delta$ and $\dhat$ as above.
We will need the following simple fact about Euclidean domains.
\begin{lem} \label{lem:delta}
  Let $a,x \in R \setminus \{0\}$ be such that
  $\delta (ax) \le \delta (a)$. Then $x$ is a unit of $R$.
\end{lem}
\begin{proof}
  There are $q, r \in R$ with
  $a = qax + r$ and $\delta(r) < \delta (ax)$.
  Then $\delta(r) < \delta(a)$. If $r = 0$, then
  $a = qax$ and thus $qx=1$ and $x$ is a unit.
  If $r \neq 0$, then 
  since $r = a (1-qx)$, we have $\delta(a) \le \delta (a (1-qx)) \le
  \delta(r)$, a contradiction.
\end{proof}  
\begin{de}[Reducibility]                     
We say that $b \in R$ is \emph{reducible} by $A \subseteq R$
if there are $a \in A$ and $q \in R$ such that
$\dhat (b - q a) < \dhat(b)$.
Now let $G \subseteq \rvxk\setminus \{0\}$. We say that a term vector
$\term{b}{\alpha}{i}$ is \emph{reducible} by $G$ if
$b$ is reducible by $\{ \LC (g) \,:\, \LM (g) \mid \monov{\alpha}{i} \}$, and
that $f \in \rvxk \setminus \{0\}$ is \emph{reducible} by $G$ if
it contains a term vector that is reducible by $G$.
\end{de}
A polynomial vector $p$ is \emph{normalized} if $p \neq 0$ and 
$\dhat (\LC(p)) \le \dhat (u \, \LC(p))$ for all units $u$ of $R$.
The subset
$G$ of $\rvxk$ is normalized if every $g \in G$ is normalized.
\begin{de} \label{de:redgb}
  Let $G$ be a strong Gr\"obner basis of the submodule $I$
  of $\rvxk$. Then $G$ is a \emph{reduced strong Gr\"obner basis}
  of $I$ if for each $g \in G$, $g$ is normalized and
  $g$ is not reducible by $G \setminus \{g\}$.
\end{de}
From an admissible ordering of the monomial vectors and the function
$\dhat$, one can define a total order on $\rvxk$. To this end,
we order polynomial vectors $p,q$ as follows: for $p \neq q$,
   let
  $\monov{\gamma}{i} := \LM (p-q)$, let
  $a$ be the coefficient of $\monov{\gamma}{i}$ in $p$,
  and let
  $b$ be the coefficient of $\monov{\gamma}{i}$ in $q$.
  Then we say $p <_P q$ if $\dhat(a) < \dhat(b)$ and
  $p \le_P q$ if $p = q$ or $p <_P q$.
  The order $\le_P$ is a well order on $\rvxk$ (Lemma~\ref{lem:lep}).
\begin{thm} \label{thm:unique}
  Let $I$ be a submodule of $\rvxk$, and
  let $\Min (\DegD(I))$ be the set of minimal elements
  of $\DegD (I)$ with respect to the ordering
  $\divd$. For every 
  $(\monoexp{\alpha}{i}, d) \in \Min (\DegD(I))$, we choose
  $g_{\alpha, i, d}$ to be the minimal element in $I$ with
  respect to $\le_P$ such that 
  $\DegD (g_{\alpha, i, d}) = (\monoexp{\alpha}{i}, d)$.
  Then 
  \[
     G := \{ g_{\alpha, i, d} \mid (\monoexp{\alpha}{i}, d) \in \Min (\DegD (I))\}
  \]
  is finite, and $G$ is the unique reduced strong Gr\"obner basis of $I$.
\end{thm}
\begin{proof}
  As an antichain in the ordered set
  $(\N_0, \le)^n \times (\{1,\ldots, k\}, {=}) \times (W, {\le})$,
  the set $\Min (\DegD(I))$ is finite
  (Theorem~\ref{thm:divdorder}\eqref{it:d2}), and hence $G$ is finite.
    As in the proof of Theorem~\ref{thm:existssgb}, we see
  that $G$ is a strong Gr\"obner basis.

  Now we show that $G$ is reduced. Let $g \in G$.
  We first show that  $g$ is normalized.
  Supposing that $g$ is not normalized, there is a unit $u \in R$ with
  $\dhat (u \, \LC(g)) < \dhat (\LC (g))$. Since
  $u$ is a unit, $\delta (u \, \LC(g)) = \delta (\LC (g))$.
  Thus $\DegD (ug) = \DegD (g)$, but $ug <_P g$. This contradicts
  the choice of $g$. Hence $g$ is normalized.

  Next, we show that $g$ is not reducible by $G \setminus \{g\}$.
  Seeking a contradiction, we suppose that
  $g$ is reducible by $G \setminus \{g\}$. Then there are a term vector
  $\term{a}{\alpha}{i}$ in $G$, $h \in G \setminus \{g\}$ and
  $q \in R$ 
    such that $\LM(h) \mid \monov{\alpha}{i}$ and 
    $\dhat (a - q \,\LC (h)) < \dhat (a)$.

 \textbf{Case 1}: $\term{a}{\alpha}{i} = \LT (g)$:
 Let $b$ be a greatest common divisor of
 $\LC(h)$ and $\LC(g)$ in $R$. Since $R$ is Euclidean,
 there exist $u,v \in R$ with $u \, \LC(h) + v \, \LC (g) = b$.
 Thus $\LT (u \frac{\LM(g)}{\LM(h)} h + v g) = \term{b}{\alpha}{i}$.
 Since $G$ is a strong Gr\"obner basis of $I$, there is
 $h_1 \in G$ with
 \begin{equation} \label{eq:h1b}
   \LT (h_1) \mid \term{b}{\alpha}{i}.
 \end{equation}  
 Since $\term{b}{\alpha}{i} \mid \LT (g)$,
 we obtain
 \begin{equation*} \label{eq:h1g}
   \LT (h_1) \mid \LT (g).
 \end{equation*}  
 Thus $\DegD (h_1) \divd \DegD (g)$.
 Since $\DegD (g) \in \Min (\DegD (G))$, we then have
 $\DegD (h_1) = \DegD (g)$.
 Since $G$ contains only one element $f$  with
 $\DegD (f) = (\monoexp{\alpha}{i}, \delta (a))$, we have
 $h_1 = g$.
 Now by~\eqref{eq:h1b}, we have $\LC (h_1) \mid b$.
 From  the definition of $b$ as a gcd, we have
 $b \mid \LC(h)$, and thus
 $\LC(h_1) \mid \LC(h)$ and therefore
 $\LC(g) \mid \LC(h)$.
 Hence there is a $q_1 \in R$ such that
 $\LC (h) = q_1 a$.
 Therefore,
 \begin{equation} \label{eq:din}
   \dhat (a- qq_1 a) < \dhat (a).
 \end{equation}
 By~\eqref{eq:dhat}, we then have
   $\delta (a- qq_1 a) \le \delta (a)$, and thus 
 by Lemma~\ref{lem:delta}, either $1-qq_1 = 0$ or
 $1-qq_1$ is a unit in $R$.
 
 \textbf{Case 1.1}: $1 - q q_1 = 0$: Then $q_1$
 is a unit in $R$ and therefore $\LC (h) \mid \LC (g)$.
 Since $\LM(h) \mid \LM(g)$, we obtain
 $\LT (h) \mid \LT (g)$ and therefore
 $\DegD (h) \divd \DegD (g)$. From the minimality of $\DegD(g)$,
 we obtain $\DegD(g) = \DegD(h)$. Since $G$ contains only one element
 $f$ with $\DegD (f) = \DegD(g)$, we have
 $h = g$, contradicting $h \in G \setminus \{g\}$.

 \textbf{Case 1.2}: $1 - q q_1$ is a unit in $R$:
 Since $g$ is normalized, we then have
 $\dhat (a (1 - qq_1)) \ge \dhat(a)$, contradicting~\eqref{eq:din}.

 \textbf{Case 2}: $\term{a}{\alpha}{i} \neq \LT (g)$:
    Then $\md (\term{a}{\alpha}{i})  < \md (g)$.
      Since $\dhat (a - q \,\LC (h)) < \dhat (a)$, we obtain
      $g - q \frac{\mono{\alpha}}{\LM(h)} h <_P g$ and
      $\LT (g - q\frac{\mono{\alpha}}{\LM(h)} h) = \LT (g)$,
      and therefore $\DegD (g - q\frac{\mono{\alpha}}{\LM(h)}) = \DegD (g)$.
      This contradicts the minimality of $g$ with respect to~$\le_P$.

      This completes the proof that $g$ is not reducible by
      $G \setminus \{g\}$.

Therefore $G$ is a reduced strong Gr\"obner basis.
The uniqueness follows from the following lemma.
\end{proof}

\begin{lem} Let $I$ be a submodule of $\rvxk$, and let $G,H$ be
  reduced strong Gr\"obner bases of $I$. Then $G = H$.
\end{lem}
\begin{proof}
  By symmetry, it is sufficient to prove $G \subseteq H$.
  Let  $g \in G$. Since $g \in I$, there is $h \in H$ such that
  $\LT (h) \mid \LT (g)$, and since $h \in I$, there
  is $g_1 \in G$ with $\LT (g_1) \mid \LT (h)$.
  If $g_1 \neq g$, then $\LC (g)$ is reducible by
  $\LC(g_1)$, contradicting the fact that $G$ is reduced.
  Thus $g_1 = g$, and therefore $\LT(g) \mid \LT(h) \mid \LT(g)$
  and thus $\LC (g)$ and $\LC(h)$ are associated in $R$.
  Since both $G$ and $H$ are normalized, we obtain $\LC(g) = \LC(h)$,
  and thus
  $\LT (g) = \LT (h)$.

  We will now show $g = h$. Seeking a contradiction, we suppose
  $g \neq h$. Since $\LT (g) = \LT (h)$, we have
  $\md (g-h) < \md (g)= \md (h)$.
  Let $\term{a}{\alpha}{i} := \LT (g - h)$, let $b$ be the coefficient of
$\monov{\alpha}{i}$ in $g$, and let $c$ be the coefficient of
$\monov{\alpha}{i}$ in $h$. Then $a = b - c$.

  Since $g-h \in I$, there is $g_1 \in G$ with $\LT (g_1) \mid \LT (g-h)$.
  We already know that $\LT(g) = \LT(h)$, and thus
  $\md (g-h) < \md (g)$.  Therefore $g_1 \neq g$.
  From $\LC (g_1) \mid b - c$, we obtain $q \in R$ such that
  $q \, \LC(g_1) = b - c$ and therefore $c = b - q \, \LC(g_1)$.
  Since $b$ is not reducible by $\{\LC(g_1)\}$, we have
  $\dhat (b) \le \dhat (c)$.  Similarly,
  since $g-h \in I$, there is $h_1 \in H$ with $\LT(h_1) \mid \LT(g-h)$.
  Since $\LT(g) = \LT(h)$, we have  $\md (g-h) < \md (h)$,
  and thus $h_1 \neq h$. From $\LC(h_1) \mid \LT(g-h)$, we obtain $q \in R$
  with $q \, \LC(h_1) = b - c$, and thus $b = c + q \, \LC (h_1)$.
  Since $c$ is not reducible by $\{\LC(h_1)\}$, we have
  $\dhat(b) \ge \dhat(c)$.
  Altogether, we have $\dhat(b) = \dhat(c)$ and thus $b = c$, leading
  to the contradiction $a=0$. Thus $g = h$ and therefore $g \in H$.
\end{proof}  

\section{Construction of reduced strong Gr\"obner bases} \label{sec:constructreduced}
We let $R$ be a Euclidean domain with grading function $\delta$
and an additional function $\dhat$ as in Section~\ref{sec:redgb}.
For reducing coefficients, we will suppose that with respect to $\dhat$, we can perform the following two algorithmic tasks:
\begin{enumerate}
\item For $a \in R$, we can find a unit $u$ in $R$ such that
  $\dhat (ua)$ is minimal in $\{\dhat(u'a) \mid u' \text{ is a unit of } R\}$.
\item For $a, b \in R$, find $q \in R$ such that
  $\dhat (b -q a)$  is minimal in $\{ \dhat (b-q'a) \mid q' \in R \}$.
\end{enumerate}
For many Euclidean domains, e.g. for the fields $\R$ or $\Q$,
it is difficult to describe such a function $\dhat$.
However, on a field $k$, it suffices to assume that
$\dhat$ satisfies $\dhat (0) < \dhat (1) < \dhat (x)$ for all
$x \in k \setminus \{0,1\}$. Then for $a \in k\setminus \{0\}$,
$u := a^{-1}$ minimizes $\dhat (ua)$ and $q := ba^{-1}$ minimizes
$\dhat (b - qa)$.
For $\Z$, we may take $\dhat (z) := 3|z| - \sgn (z)$, which
yields
$\dhat(0) < \dhat (1) < \dhat (-1) < \dhat (2) < \dhat (-2) < \cdots$.
Then for $a \in \Z \setminus \{0\}$, $u := \sgn (a)$ minimizes
$\dhat (ua)$ and the unique $q \in \Z$ with
$-\frac{a}{2} < b - q a \le \frac{a}{2}$ minimizes
$\dhat (b - qa)$.

For reducing polynomial vectors, we follow \cite{Li:ECOS}
and use reductions that, when they affect
the leading term of polynomial, eliminate this leading term
in one step. We call such reductions \emph{soft}.
\begin{de} \label{de:softred}
  The polynomial vector
  $f \in \rvxk \setminus \{0\}$ is \emph{softly reducible} by $G$ if
  $f - \LT (f)$ is reducible by $G$ or there is $g \in G$ such
  that $\LT (g) \mid \LT(f)$.
  A set $G \subseteq \rvxk$ is \emph{softly reduced} if no
  $f \in G$ is softly reducible by $G \setminus \{f\}$.
\end{de}
We will consider the following ordering of finite subsets
of $\rvxk$. We say that 
  $G_1 \le_S G_2$  if there is an injective map $\phi : G_1 \to G_2$
such that $g \le_P \phi(g)$ for all $g \in G_1$.
This ordering is a well partial ordering (Lemma~\ref{lem:les}).
One step of a soft reduction is performed in the following
algorithm \textsc{SoftlyReduce}.
\begin{alg}[Soft reduction] \label{alg:sr} \mbox{}\\
  Input: $G \subseteq \rvxk \setminus \{0\}$
          such that $G$ is not softly reduced. \\
  Output: $H \subseteq \rvxk \setminus \{0\}$ 
  such that
  \begin{enumerate}
  \item $\submod{H} = \submod{G}$,
  \item Every $g \in G$ has a strong standard representation by $H$,
  \item $H <_S G$
  \end{enumerate}  
  \begin{algorithmic}[1]
    \Function{SoftlyReduce}{$G$}
    \State Choose $f, h \in F$ and 
       a term $\term{a}{\alpha}{i}$ from $f$ such that $f \neq h$, 
       $\LM(h) \mid \monov{\alpha}{i}$ and
       there is $q \in R$ such that
       ($\monov{\alpha}{i} = \LM (f)$ and $a - q \, \LC (h) = 0$)
        or 
        ($\monov{\alpha}{i} \neq \LM (f)$ and
        $\dhat(a - q \, \LC (h)) < \dhat (a)$).
     \State $r \gets f - q \frac{\monov{\alpha}{i}}{\LM(h)} h$
     \If{$r = 0$}
     \State $H \gets G \setminus \{f\}$
     \Else
      \State $H \gets (G \setminus \{f\}) \cup \{r\}$ 
     \EndIf  
     \State Return $H$
     \EndFunction
   \end{algorithmic}      
\end{alg}
\begin{lem}
  Algorithm~\ref{alg:sr} is correct.
\end{lem}  
\begin{proof}
  Clearly, $H$ and $G$ generate the same submodule.

  Next, we show that every $g \in G$ has a strong standard representation.
  Let $g \in G$. If $g \in H$, then
  $g = 1 \mono{0} g$ is such a representation.
  If $g \not\in H$, then $g = f$ and
  \begin{equation} \label{eq:fqr}
    f = q \frac{\monov{\alpha}{i}}{\LM(h)} h + 1 \mono{0} r.
  \end{equation}
  We first assume $r \neq 0$.
  If $\LT (f) = \term{a}{\alpha}{i}$, then
  $\md (\frac{\monov{\alpha}{i}}{\LM(h)} h) = (\alpha, i) = \md (f)$
  and
  $\md (\mono{0} r) = \md (r) < \md (f)$.
  If $\LT (f) \neq \term{a}{\alpha}{i}$, then
  $\md (\frac{\monov{\alpha}{i}}{\LM(h)} h) = (\alpha, i) < \md (f)$ and
  $\md (\mono{0} r) = \md (r) = \md (f)$.
  In both cases~\eqref{eq:fqr} is a strong standard representation
  of $f$ by $H$ with remainder~$0$.
  If $r = 0$, then  $f = q \frac{\monov{\alpha}{i}}{\LM(h)} h$
  is a strong standard representation.

  For proving $H <_S G$, we define $\phi : H \to G$
  by $\phi (h) = h$ for $h \in H \setminus \{r\}$,
  and $\phi (r) = f$ when $r \neq 0$.
  Since $r <_P f$, the mapping $\phi$ witnesses $H <_S G$.
\end{proof}
We also need to normalize polynomial vectors. The following
procedure normalizes one vector in $G$.
\begin{alg}[Normalization] \label{alg:norm} \mbox{}\\
  Input: $G \subseteq \rvxk \setminus \{0\}$ such that
         $G$ contains an element that is not normalized.
  Output: $H \subseteq \rvxk \setminus \{0\}, H \neq \emptyset$ 
  such that
  \begin{enumerate}
  \item $\submod{H} = \submod{G}$,
  \item Every $g \in G$ has a strong standard representation by $H$,
      \item $H <_P G$.
  \end{enumerate}
  \begin{algorithmic}[1]
    \Function{Normalize}{$G$}
    \State Choose $g \in G$ such that $g$ is  not normalized
    \State Find a unit $u$ in $R$ such that $ug$ is normalized
    \State $H \gets (G \setminus \{g\}) \cup \{ug\}$
    \State Return $H$
  \EndFunction   
  \end{algorithmic}
\end{alg}

\begin{lem} Algorithm~\ref{alg:norm} is correct.
\end{lem}
\begin{proof}
  It is clear that $H$ and $G$ generate the same submodule.

  Furthermore, $g = u^{-1} \mono{0} (ug)$ is a strong standard representation
  of $g$ by $H$.

  We have $ug <_P g$. Hence
  $\phi (h) := h$ for $h \in H \setminus \{ug\}$ and
  $\phi (ug) = g$ witnesses $H <_S G$.
\end{proof}  

\begin{thm} \label{thm:sred} Let $G$ be a softly reduced strong Gr\"obner basis in
  which every element is normalized. Then $G$ is reduced.
\end{thm}
\begin{proof}
  Let $g \in G$. We have to show that $g$ is not reducible
  by $G \setminus \{g\}$. Suppose that $g$ is reducible. Then
  there are $h \in G \setminus \{g\}$, a term
  $\term{a}{\alpha}{i}$ in $g$ and $q \in R$ such that
  $\LM (h) \mid \monov{\alpha}{i}$ and $\dhat (a - q \, \LC (h))
  < \dhat (a)$.
  If $\term{a}{\alpha}{i} \neq \LT (g)$, then $g - \LT (g)$ is
  reducible by $\{h\}$, and thus $g$ is softly reducible by
  $G \setminus \{g\}$.
  If $\term{a}{\alpha}{i} = \LT (g)$, then
  let $d := \gcd (\LC(g), \LC(h))$. 
  There is a polynomial vector $f$ in the module generated by $G$
  such that $\LT (f) = \term{d}{\alpha}{i}$ and thus there
  is $g_1 \in G$ with $\LT (g_1) \mid \term{d}{\alpha}{i}$.
  Since $G$ is softly reduced, we then have $g_1 = g$, and thus
  $a = \LC (g) = \LC (g_1) \mid  \LC (h)$.
  Then $\LC (g) - q \, \LC(h)$ is a multiple of~$a$.
  Since $\dhat (a - q \, \LC (h)) < \dhat (a)$,~\eqref{eq:dhat}
  implies
  $\delta (a-q\,\LC(h)) \le \delta (a)$, and
  thus 
  by Lemma~\ref{lem:delta}, there is a unit in $R$
  such that $u \, \LC(g) = \LC(g) - q \, \LC (h)$.
  Since $g$ is normalized, $\dhat (\LC (g)) \le \dhat (\LC(g) - q \, \LC (h))$;
  this contradicts $\dhat (\LC(g) - q \, \LC (h)) < \dhat (\LC (g))$.
\end{proof}  
  
In the computation of a strong Gr\"obner basis, we may interleave
the three steps done in \textsc{Augment}, \textsc{SoftlyReduce}
and \textsc{Normalize} as we wish. However, at some point,
we may for instance enter the while loop with $G$ normalized and softly reduced:
then in this course of the while-loop, we have to use
the procedure \textsc{Augment}. Note that the while-condition
guarantees that we have at least one choice in every execution
of the while-loop.
\begin{alg}[Reduced Strong Gr\"obner Basis] \label{alg:rsgb}
       \mbox{} \\
       Input: $F \subseteq \rvxk \setminus \{0\}$,
         an admissible order $\le$ of $\NONk$. \\
       Output: $G \subseteq \rvxk \setminus \{0\}$ such that
       $G$ is a reduced strong Gr\"obner basis of the submodule
       generated by $F$ with respect to the monomial vector ordering
       $\le$.
        \begin{algorithmic}
         \State $G \gets F$
         \State $P \gets \emptyset$
         \While{($\exists p, q \in G : p \neq q$  and $\{p, q\} \not\in P$)
           or  \\ \mbox{\,\,\,\,\,\,\,\,\,\,\,\,\,\,\,\,}  ($G$ is not softly reduced) or \\
                  \mbox{\,\,\,\,\,\,\,\,\,\,\,\,\,\,\,\,} ($G$ is not normalized)}
         \State \textbf{Do exactly one} out of the possible choices from
                (1),(2),(3):
         \State (1) $(G, x) \gets \textsc{Augment} (G, \{p,q\})$
         \State \mbox{\,\,\,\,\,\,\,\,\,} \textbf{if} {$x = 1$} \textbf{then}  $P \gets P \cup \{\{p,q\}\}$ %
         \State (2)  $G \gets \textsc{SoftlyReduce} (G)$
         \State (3)  $G \gets \textsc{Normalize} (G)$
      \EndWhile    
       \State Return $G$
   \end{algorithmic}
\end{alg}

\begin{thm} Algorithm~\ref{alg:rsgb} terminates on every input and
  produces a correct result.
\end{thm}

\begin{proof}
  We first show that the algorithm terminates.
  Seeking a contradiction, we consider an execution that
  runs forever. In this execution,
  let $G_i$ be the value of $G$ at the beginning of
  the $i$th execution of the while-loop. The output conditions of
  the three algorithms \textsc{Augment}, \textsc{SoftlyReduce}
  and \textsc{Normalize} imply that
  $\ups{\DegD(G_i)} \subseteq \ups{\DegD (G_{i+1})}$.
  Thus there is an $n_1 \in \N$ such that for all $i \ge n_1$,
  we have 
  $\ups{\DegD (G_{i})} = \ups{\DegD (G_{i+1})}$.

  From this point onwards, the assignments to $G$ in
  lines~\ref{li:assg1},~\ref{li:assg2},~\ref{li:assg3} in \textsc{Augment}
   (Algorithm~\ref{alg:augment}) will not
  be executed any more because all of these assignments strictly
  increase $\ups{\DegD(G)}$ with respect to $\subseteq$.
  In other words, $G$ will not be changed any more by
  \textsc{Augment}, which also follows from output condition~\eqref{it:aug2}
  of \textsc{Augment}.
  Hence for all $i \ge n_1$, we have $G_{i+1} \le_S G_i$.
  Thus there is $n_2 \in \N$ with $n_2 \ge n_1$ such that for all $i \ge n_2$,
  $G_{i+1} = G_i$.
  From this point on, \textsc{SoftlyReduce} and
  \textsc{Normalize} cannot be called any more because
  both of them strictly decrease $G$ with respect to $\le_S$.
  Hence, the only remaining possible branches are
  the cases $f = 0$ and and
  $\exists g \in G : \LT(f) \mid \LT (g)$ in the execution
  of \textsc{Augment}.  In detail, only
  the assignments contained in line~\ref{li:e1} to~\ref{li:e2}
  of \textsc{Augment}
  can be excuted. In both branches $x = 1$ (this can also be seen directly from
  output condition~\eqref{it:aug3} of \textsc{Augment}), and thus
  $\# ({G_{i+1} \choose 2} \setminus P_{i+1}) < \# ({G_i \choose 2} \setminus P_i)$.  Hence, starting from the $n_2$th execution of the while-loop
  of Algorithm~\ref{alg:rsgb}, this nonnegative number strictly decreases
  forever, which is impossible.
  Hence the algorithm terminates on every input.

  From Lemma~\ref{lem:transssr}, we obtain
  that throughout the execution of the algorithm, the
  set
  $\{ f \in \rvxk \mid f$ has a strong standard representation
  by $G \}$ increases with respect to $\subseteq$.
  By the output conditions of all three procedures
  \textsc{Augment}, \textsc{SoftlyReduce} and \textsc{Normalize},
  $\submod{G} = \submod{F}$.
  Therefore, when the while-loop is left, $G$ is softly reduced
  and $G$ is normalized. Furthermore, every two-element subset
  $\{p, q\}$ of $G$ lies in $P$ and therefore
  has an $S$-polynomial vector that is $0$ or has a strong standard
  representation by $G$.
  Thus by Theorem~\ref{thm:Spols}, $G$ is a strong Gr\"obner basis
  of $\submod{G}$, and by Theorem~\ref{thm:sred}, $G$ is
  reduced.
\end{proof}

\section{Linear algebra over $\rvx$} \label{sec:la}
Let $D$ be a commutative ring with unit. By $D^{r \times s}$,
we denote the set of $r \times s$-matrices over $D$.
For $A \in D^{r \times s}$, we define 
$\col(A) = \{ A x \mid x \in D^{s} \}$ as the  \emph{column module}
and  $\row (A) = \{ y A \mid y \in D^{r} \}$ as
the \emph{row module} of  $A$.
The set $\ker (A) = \{ y \in D^{s} \mid  A y = 0 \}$
is the \emph{kernel} oder \emph{null module} of $A$.
We will now compute bases for these modules
in the case $D = \rxn$, where $R$ is a Euclidean domain.
We assume that we have the Euclidean grading function $\delta$
and $\dhat$ for $R$ as in Section~\ref{sec:redgb}.
As an additional assumption, we assume
that $\dhat(1)$ is minimal in $\{\dhat(u) \mid u$ is a unit of
$R\}$.
For a matrix $A \in \rxn^{r \times s}$ and admissible monomial orders
$\le_1, \ldots, \le_s$ on the monomials of $\rvx$, we define the
\emph{position over term}-order $\le$ by
$\monov{\alpha}{i} \le \monov{\beta}{j}$ if
$i > j$ or $(i = j \text{ and } \alpha \le_i \beta)$.
We say that a matrix $H \in \rvx^{r \times s}$ is the \emph{Gr\"obner normal form} with respect to $(\le_1, \ldots, \le_s)$
for $A$ 
if the rows of $H$ are a reduced strong Gr\"obner basis of
the module $\row (A)$ with respect to $\le$, and the rows are ordered
in strictly decreasing order with respect to the total order $\le_P$
defined after Definition~\ref{de:redgb}. An example
of such a matrix is given in \eqref{eq:H}.
The entries of $H$ can be described as follows:
\begin{lem} \label{lem:fork}
    Let $R$ be a Euclidean domain, 
    let $A \in \rvx^{r' \times s}$,
    and let $H = (h_{ij})_{(i,j) \in \ul{r} \times \ul{s}}\in \rvx^{r \times s}$
    be the Gr\"obner normal form of $A$.
    For $i \in \ul{s}$, we define the
    \emph{$i$th step} of $H$ by
    \[
    S_i = \{ h_{t,i} \setsuchthat t \in \ul{r},  h_{t,i} \neq 0, \text{ and }
    h_{t, 1} = \cdots = h_{t, i-1}  = 0 \}.
    \]
    The \emph{$i$th fork ideal} of $\row (A)$ is the set
    \[
    F_i = \{ p \in \rvx \mid \exists p_{i+1}, \ldots, p_s \in \kvx :
    (\underbrace{0,\ldots, 0}_{i-1}, \, p, \,  p_{i+1}, \ldots, p_s)
    \in \row (A) \}.
    \]
    Then $S_i$ is a reduced strong Gr\"obner basis of the ideal $F_i$ of $\rvx$
    with 
    respect {to~$\le_i$}.
 \end{lem}
 \begin{proof}
   Let $p \in F_i$ with $p \neq 0$, and let $\vb{v} = (\underbrace{0,\ldots, 0}_{i-1}, \, p, \, p_{i+1}, \ldots, p_s) \in \row(A)$.
   Then $\vb{v} = p \,e_i + \sum_{j=i+1}^s p_j \, e_j$.
   Let $h_1, \ldots, h_r$ be the rows of $H$.
   Since $\{h_1, \ldots, h_r\}$ is a strong Gr\"obner basis
   of $\row(A)$, there is %
   $t \in \ul{r}$ such that  $\LT (h_t) \mid \LT (\vb{v}) = \LT(p \, e_i)$.
   Then $\LT(h_t)$ is of the form
   $\term{a}{\alpha}{i}$, and therefore 
  $h_{t,i} \in S_i$. Hence  $h_t$ can be written as
  $(0,\ldots, 0, h_{t,i}, h_{t, i+1}, \ldots, h_{t,s})$ with
   $\LT (h_t) = \lt(h_{t,i}) e_i$. (Recall from Section~\ref{sec:construct}
   that we write
   $\LT(f)$ when $f$ is a polynomial vector in $\rvxk$ and $\lt(f)$ when
   $f$ is a single polynomial in $\rvx$.)
  Hence $\lt(h_{t,i}) \mid \lt (p)$. Thus $S_i$
  is a Gr\"obner basis of $F_i$.

  Now suppose that $S_i$ is not reduced. Then we have
  $h_{u,i}, h_{v,i} \in S_i$ with $u \neq v$, $q \in R$ and
  $\alpha \in \NON$ such that
  $\deg (h_{u, i}) \ge \deg (h_{v, i})$ and
  $h_{u, i} >_p h_{u,i} - q \mono{\alpha} h_{v, i}$, where $\le_p$ is
  defined from $\le_i$ for polynomials in analogy to the definition
  of $\le_P$ for polynomial vectors in Section~\ref{sec:redgb}.
    Then $h_u >_P h_u - q \mono{\alpha} h_v$, contradicting the fact that
    the rows of $H$ are a reduced Gr\"obner basis.
 \end{proof}  
 This allows us to solve  linear systems over
 $\rvx$.
 As an example, we consider the linear equation
 $(10 y) z_1 + 0 z_2 + (4 x) z_3 = 4 x^3$, where we
 look for the set of all solutions $(z_1, z_2, z_3) \in \Z[x,y]^3$.
 We collect the data from this equation in the matrix
 \[
     A' := \left(
\begin{array}{ccccc}
 -4 x^3 & 1 & 0 & 0 & 0 \\
 10 y & 0 & 1 & 0 & 0 \\
 0 & 0 & 0 & 1 & 0 \\
 4 x & 0 & 0 & 0 & 1 \\
\end{array}
\right)
\]
and we compute the Gr\"obner normal form (with respect to the lexicographical
ordering with $x > y$ in all columns) of $A'$ as
\begin{equation} \label{eq:H}
  H =
\left(
\begin{array}{ccccc}
 2 x y & 0 & x & 0 & -2 y \\
 4 x & 0 & 0 & 0 & 1 \\
 10 y & 0 & 1 & 0 & 0 \\
 0 & 1 & 0 & 0 & x^2 \\
 0 & 0 & 2 x & 0 & -5 y \\
 0 & 0 & 0 & 1 & 0 \\
\end{array}
\right).
\end{equation}
Then we can read from this matrix that $(0,0,x^2)$ is one solution,
and the solution module of  $(10 y) z_1 + 0 z_2 + (4 x) z_3 = 0$
is generated, as a $\Z[x,y]$-module,
by $(2 x, 0,-5 y)$ and $(0, 1, 0)$.
This is justified by the following theorem, which explains how to
solve linear systems in a style that follows \cite[Chapter~3.8]{AL:AITG}.
 \begin{thm} \label{sat:linsys}
    Let $R$ be a Euclidean domain, 
    let $A \in \rvx^{r \times s}$, let $b \in \rvx^{r \times 1}$ and
    let $\le_{-1}, \le_0, \le_1, \ldots, \le_s$ be admissible orders
    on the monomials of $\rvx$.
    Let $H \in \rvx^{r' \times (r + s + 1)}$ be the Gr\"obner normal form of
    \[
     A' = \left(
        \begin{array}{c|c}
          \begin{array}{c}
            -b^T  \\
            A^T
          \end{array} &
             I_{s + 1}
        \end{array}
        \right)
        \]
        with respect to the monomial orders $(\underbrace{\le_{-1}, \ldots, \le_{-1}}_{r \text{ times}},
        \le_0, \le_1, \ldots, \le_s)$.
    We write $H$ as
    \[
    H =
    \left(
    \begin{array}{ccc}
      B & * & * \\
      0 & v & S \\
      0 & 0 & D
    \end{array}
    \right),
    \]
    where $B$ has exactly $r$ columns, $v$ exactly $1$ column and 
    $D$ exactly  $s$ columns, and furthermore
    the last line of $B$ is not the zero-vector,
    and the last entry of $v$ is not~$0$.
    Then we have:
    \begin{enumerate}
       \item \label{it:v}
         The entries of $v$ are a reduced strong Gr\"obner basis
         of the ideal
         \[
         (\col(A) : b) := \{ p \in \rvx \mid p \, b \in \col (A) \}.
         \]
         of $\rvx$ with respect to the monomial order $\le_0$.
       \item  \label {it:S}  
         The system $Ax = b$ has a solution
         in $\rvx^{s}$ if and only if  $v=(1)$.
         Then the matrix $S$ has exactly one row $s_1$, and $s_1$ is
         the minimal solution of $Ax = b$ with respect to $\le_P$, where
         $\le_P$ is the total order on polynomial vectors defined from
         the admissible order $\le$ that is the position over term order
         coming from $(\le_1, \ldots, \le_s)$.
      \item  \label{it:D}
        $D$ is in Gr\"obner normal form and 
        $\row (D) = \ker (A)$.
    \end{enumerate}
\end{thm}    
\begin{proof}
  \eqref{it:v}
  We first show that $\{p \in \rvx \mid p\,b \in \col(A)\}$ is
  equal to the $(r+1)$th fork ideal $F_{r+1}$ of $A'$.
  To this end, let $a_1, \ldots, a_s \in \rvx^r$ be the column
  vectors of $A$.
  For proving one inclusion, we assume that $p_{r+1} \in F_{r+1}$. Then there
  are $p_{r+2}, \ldots, p_{r+s+1} \in \rvx$ such that
  $(0, \ldots, 0, p_{r+1}, p_{r+2}, \ldots, p_{r+s+1})$ is in $\row(A')$,
  and thus there is $(f_0,f_1, \ldots, f_s) \in \rvx^{s+1}$
  such that
  \begin{equation} \label{eq:aprime}
  (f_0,f_1, \ldots, f_s) \cdot A' = (0, \ldots, 0, p_{r+1}, p_{r+2}, \ldots, p_{r+s+1}).
  \end{equation}
  Considering the first $r$ entries of the right hand side
  of~\eqref{eq:aprime}, we obtain $-f_0 b + \sum_{i=1}^s f_i a_i = 0$,
  and hence  $f_0 b = \sum_{i=1}^s f_i a_i$, and
  therefore $f_0 b \in \col(A)$.
  The $(r+1)$th column of $A'$ is the first unit vector in
  $\rvx^{s+1}$. Hence $f_0 = p_{r+1}$, und thus $p_{r+1} b \in \col (A)$
  and therefore $p_{r+1} \in (\col (A) : b)$.

  Now assume that $p \in (\col (A) : b)$. Then there is
  $(f_1, \ldots, f_s) \in \kvx^s$ such that
  $\sum_{i=1}^s f_i a_i = p b$. Therefore the first
  $r$ columns of 
  $(p, f_1, \ldots, f_s) \cdot A'$ are $0$, and therefore
  the $(r+1)$th entry of $(p, f_1, \ldots, f_s) \cdot A'$
  is an element of $F_{r+1}$. Since this entry is $p$, we have
  $p \in F_{r+1}$.

  By Lemma~\ref{lem:fork}, the entries of $v$ are a reduced
  strong Gr\"obner basis of $F_{r+1} = (\col (A) : b)$ with respect to $\le_0$.
  
\eqref{it:S}
The system $Ax=b$ has a solution if and only if
$b \in \col (A)$, which means $1 \in (\col(A) : b)$.
This holds if and only if the reduced Gr\"obner basis of
$(\col (A) : b)$ is $\{1\}$.
By item~\eqref{it:v}, the entries of $v$ are a reduced Gr\"obner
basis of $(\col (A) : b)$.
Altogether, $Ax= b$ has a solution in $\rvx^s$ if and only if $v = (1)$.

\eqref{it:D}
It is not hard to show that the rows of $D$ are a reduced strong
Gr\"obner basis of the module
\[
E :=\{(f_{r+2}, \ldots, f_{r+s+1}) \in \rvx^{s}
\mid (0,\ldots,0, f_{r+2}, \ldots, f_{r+s+1}) \in \row (A')\}.
\]
We now show $E = \ker (A)$.
For $\subseteq$, we assume $(0,\ldots, 0, f_{r+2}, \ldots, f_{r+s+1}) \in \row(A')$.
Then there is
$(g_0, g_1, \ldots, g_s) \in \rvx^{s+1}$ with
\begin{equation} \label{eq:gA}
  (g_0, g_1, \ldots, g_s) \cdot A' = (0, \ldots, 0, f_{r+2}, \ldots, f_{r+s+1}).
\end{equation}  
Hence $g_0 = 0$ and $(g_1, \ldots, g_s) \cdot A^T = 0$,
and therefore $(g_1, \ldots, g_s) \in \ker (A)$. Since
$(g_1, \ldots, g_s) = (f_{r+2}, \ldots, f_{r+s+1})$,
we obtain that $(f_{r+2}, \ldots, f_{r+s+1}) \in \ker (A)$.

If $(g_1, \ldots, g_s) \in \ker (A)$, then
$(0,g_1, \ldots, g_s) \cdot A' = (0,\ldots, 0, g_1, \ldots, g_s)$
and thus $(g_1, \ldots, g_s) \in E$.
\end{proof}

The Gr\"obner normal form generalizes the row echelon normal form of a matrix  $A$
over a field $k$ as computed, e.g., in Mathematica \cite{WR:MV14} by {\tt RowReduce [A]}.
To see this, we set
$R := k$ and  consider $A$ as a matrix over $R[x_1]$ (in which $x_1$ never appears).
Similarly, it also generalizes the Hermite normal form of a matrix over $\Z$ (with the
elements above the pivot elements normalized to minimize their absolute values, and preferring
$3$ over $-3$).
Here we consider $A$ as a matrix over $\Z[x_1]$, and set  $R := \Z$,
$\delta (z) := |z|$ and $\dhat (z) = 3 |z| - \sgn(z)$ for all $z \in \Z$ to obtain
$\dhat (0) < \dhat (1) < \dhat (-1) < \dhat (2) < \dhat (-2) < \cdots$.  
Hence Theorem~\ref{thm:unique} also implies the uniqueness of these normal forms.

\section{Partial orders} \label{sec:order}

A partially ordered set $(A, \rho)$ is a set $A$ together with
a partial order, i.e.,
a reflexive, transitive and antisymmetric relation $\rho$.
Often, we write $a \le b$ or $b \ge a$ for $(a, b) \in \rho$,
and $a < b$ or $b > a$ when $(a, b) \in \rho$ and $(b,a) \not\in \rho$.
We say that $a$ and $b$ are \emph{uncomparable} and write
$a \perp b$ if $(a, b) \not\in \rho$
and $(b,a) \not \in \rho$.
The sequence $(a_i)_{i \in \N}$ is an \emph{infinite descending chain}
in $A$ when $a_{i} > a_{i+1}$ for all $i \in \N$, and an
\emph{infinite antichain} when $a_i \perp a_j$ for all
$i, j \in \N$ with $i \neq j$.  An order relation $\le$ on $A$
is a \emph{well partial order} if it has no infinite descending
chains and no infinite antichains. It is a \emph{well order}
if it is total (i.e., has no distinct uncomparable elements)
and has no infinite descending chains. For a subset $B$ of
the partially ordered set $(A, \le)$, $b \in B$ is \emph{minimal}
in $B$ if there is no $b'\in B$ with $b' < b$. The subset $B$
is \emph{upward closed} if for all $b \in B$ and $a \in A$ with
$b \le a$, we have $a \in B$.
The product of $(A_1, \rho_1)$ and $(A_2, \rho_2)$ is the
set $A_1 \times A_2$ ordered by the relation
$\rho$ defined by $((a_1, a_2), (b_1, b_2)) \in \rho
:\Leftrightarrow (a_1, b_1) \in \rho_1$ and $(a_2, b_2) \in \rho_2$.
Our investigation of these orderings is facilitated by
Ramsey's Theorem \cite{Ra:OAPO} (cf. \cite{Ne:RT}):
Denote the two element subsets of $\N$ by ${\N \choose 2}$ and
let $c$ be a function from ${\N \choose 2}$ into a finite set.
Then there exists an infinite subset $T$ of $\N$ such that
$c$ is constant on ${T \choose 2}$.
All results in this section are well known; some
are taken from the survey~\cite{AA:DLHT}.
\begin{thm}[Dickson's Lemma \cite{Di:FOTO}] \label{thm:dick}
  The product of two well partially ordered sets is well partially
  ordered.
\end{thm}
\begin{proof}
  Let $(A, \le_A)$ and $(B, \le_B)$ be well partially ordered sets,
  and let $((a_i, b_i))_{i \in \N}$ be any sequence from $A \times B$.
  We colour the two element subsets of $\N$ with one of the
  nine colours from $\{\le, >, \perp\}^2$ as follows:
  when $i < j$ then
  $C (\{i,j\}) = (\le, \le)$ if $a_i \le a_j$ and $b_i \le b_j$,
  $C (\{i,j\}) = (\le, >)$  if $a_i \le a_j$ and $b_i > b_j$, \ldots .
  By Ramsey's Theorem
  there is an infinite subset $T$
  of $\N$ such that all two-element subsets of $T$ have the
  same color $c$. If this colour $c$ is not $(\le, \le)$, then we
  find an infinite descending chain or an infinite antichain
  in either $A$ or $B$. Hence $c = (\le, \le)$.
  This implies that $((a_i, b_i))_{i \in \N}$ is neither an infinite
  descending chain nor an infinite antichain.
\end{proof}
\begin{lem} \label{lem:adwell}
  Let $\le_{\operatorname{a}}$ be an admissible ordering on $\monosnk$. Then
  there is no infinite descending chain
  $m_1 >_{\operatorname{a}} m_2 >_{\operatorname{a}} \cdots$  with
  respect to this ordering.
\end{lem}
\begin{proof}
Let $(m_i)_{i \in \N}$ be a sequence from $\monosnk$.
We colour two-element subsets $\{i,j\}$ of $\N$ with $i < j$
by $C (\{i,j\}) = 1$ if $m_i \mid m_j$, $C(\{i,j\}) = 2$ if
$m_j \mid m_i$ and $m_j \neq m_i$, and $C(\{i,j\}) = 3$ if $m_i \nmid m_j$ and
$m_j \nmid m_j$.
We use Ramsey's Theorem to obtain an infinite subset $T$ of $\N$
such that all two-element subsets of $T$ have the same colour $c$.
If this colour is $2$ or $3$, then we obtain an infinite descending
chain or an infinite antichain in $\monosnk$, which is order isomorphic
to $(\N_0, {\le})^n \times (\{1,\ldots, k \},{=})$, contradicting
Theorem~\ref{thm:dick}. Hence this colour is $1$ and thus
there are $i, j \in \N$ with $i < j$ such that
$m_i \mid m_j$. Then $m_i \le_{\operatorname{a}} m_j$.
Thus $(m_i)_{i \in \N}$ cannot be an infinite descending chain.
\end{proof}

As another  consequence, we obtain that
the order relation $\divd$ defined in~\eqref{eq:lmd},
which is the order of the direct product
of $n$ copies of  $(\N_0,{\le})$ with
$(\{1,\ldots, k\},{=})$ and $(W,{\le})$ has no infinite descending
chain and no infinite antichain:
\begin{thm} \label{thm:divdorder}
   Let $(W, \le)$ be a well ordered set, and  
   let $\divd$ be the ordering on $\NONk \times W$ defined
   in~\eqref{eq:lmd}. Then we have
   \begin{enumerate}
   \item \label{it:d1} The order $\divd$ is a well partial order.
   \item \label{it:d2} For every subset $D$ of $\NONk \times W$, the set
     $\Min (D)$ of minimal elements of $D$ with respect
     to $\divd$ is finite, and for every $d \in D$ there
     is $d' \in \Min (D)$ with $d' \divd d$.
   \item \label{it:d3}
     There is no infinite ascending chain
     $D_1 \subset D_2 \subset \cdots$ of upward closed
     subsets of $(\NONk \times W, \divd)$.
\end{enumerate}     
\end{thm}
\begin{proof}
  \eqref{it:d1}
  $(\NONk \times W, \divd)$ is order isomorphic to the
  product of $n$ copies of $(\N_0,{\le})$ with
  $(\ul{k}, =)$ and $(W, \le)$. Since all factors
  are well partially ordered, Theorem~\ref{thm:dick}
  implies that $\divd$ is a well partial order.
  \eqref{it:d2}
  Distinct minimal elements of $D$ of are all uncomparable
  with respect $\divd$. Since
  $\divd$ is a well partial order and therefore has
  no infinite antichains, $\Min (D)$ is finite.
  Now let $d \in D$. If $\{ x \in D \mid x \divd d\}$
  has no minimal element, we can construct a sequence
  $(d_i)_{i \in \N}$ with
  $d \sqsupset_{\delta} d_1 \sqsupset_{\delta} d_2 \sqsupset_{\delta} \cdots$
  of elements from $D$; such sequences do not exist because
  $\divd$ is a well partial order, and therefore
  $\{x \in D \mid x \divd d\}$ has a minimal element, which
  is then also minimal in $D$.
  \eqref{it:d3}
  Let $U_1 \subset U_2 \subset \cdots$ be an infinite ascending
  chain of upward closed subsets of $\NONk \times W$.
  Then $U := \bigcup_{i \in \N} U_i$ has a finite set of minimal
  elements $\Min (U)$. Thus there is $k \in \N$ with
  $\Min (U) \subseteq U_k$, and therefore $U \subseteq U_k$,
  which yields the contradiction $U_{k+1} \subseteq U_k$.
\end{proof}  

Next, we see that the order $\le_P$ of polynomial vectors
defined before Theorem~\ref{thm:unique} is a well order.
For $f = \sum_{(\alpha, i) \in E} \term{c_{(\alpha, i)}}{\alpha}{i}$,
we let $[\monov{\alpha}{i}] \, f := c_{(\alpha, i)}$ denote the
coefficient of $\mono{\alpha}$ of the $i$\,th component of $f$.
Then for $p \neq q$, we have
$p <_P q$ if $\dhat ([\LM(p-q)]\,p) < \dhat ([\LM(p-q)]\,q)$ and
  $p \le_P q$ if $p = q$ or $p <_P q$. 
\begin{lem} \label{lem:lep}
  The relation $\le_P$ is a well order on $\rvxk$.
\end{lem}
\begin{proof}
  The relation $\le_P$ is clearly reflexive and
  antisymmetric. For transitivity, we assume
  $p <_P q <_P r$. Then
  $\LM (p - r) = \LM ((p - q) + (q - r))$, and thus
  $\md (p-r) \le \max (\md (p-q), \md (q-r))$.
  If $\LM (p-q) = \LM (q-r)$, then
  $\dhat ([\LM (p-q)] \, p) < \dhat ([\LM (p-q)] \, q) <
  \dhat ([\LM (p-q)] \, r)$.
  Hence $[\LM (p-q)] \, (p - r) \neq 0$. Thus
  $\md (p - r) \ge \md (p-q)$, and therefore
  $\md (p-r) = \md (p - q)$. Now
  $\dhat ([\LM (p-r)] \, p) < \dhat ([\LM (p-r)] \, r)$, and thus
    $p <_P r$.
    If $\LM (p-q) \neq \LM (q-r)$, then we first consider
    the case $\md (p-q) < \md (q-r)$. Then
    $\md (p-r) = \md ((p-q) + (q-r)) = \md (q-r)$.
    Since $\md(p-q) < \md (q-r)$, we have
    $[\LM(q-r)]\, p = [\LM(q-r)] \, q$. Therefore
    $\dhat ([\LM(q-r)]\, p) < \dhat([\LM(q-r)] \, r)$ and thus
    $p <_P r$. The case $\md (p-q) > \md (q-r)$ is similar.
    Thus $\le_P$ is transitive.
    It is easy to see that the ordering $\le_P$ is total.
    Now let $(f_i)_{i \in \N}$ be an infinite descending
    chain with respect to $\le_P$; among such chains, choose
    one for which $\md (f_1)$ is minimal.
    Then we must have $\md (f_1) = \md (f_i)$ for all $i \in \N$,
    since otherwise $(f_{j})_{j \ge i}$ would contradict
    the minimality. Thus $\dhat (\LC(f_i)) \ge \dhat (\LC(f_{i+1}))$
    for all $i \in \N$. This means that there is $n_1 \in \N$
    such that $\dhat (\LC(f_i)) = \dhat (\LC(f_{i+1}))$ and therefore
    $\LT (f_i) = \LT (f_{i+1})$ for all
    $i \ge n_1$. Then $(f_i - \LT (f_i))_{i \ge n_1}$ is an infinite
    descending sequence with respect to $\le_P$, contradicting the minimality
    of $\md (f_1)$.
 \end{proof}    
After stating Definition~\ref{de:softred},
we have ordered finite subsets $F,G$ of $\rvxk$ by
$F \le_S G$ if there is an injective $\phi : F \to G$ with
$f \le_P \phi (f)$ for all $f \in F$.
\begin{lem} \label{lem:les}
  The relation $\le_S$ is a well order on $\rvxk$.
\end{lem}  
\begin{proof}
  Reflexivity and transitivity of $\le_S$ are immediate.
  For checking that $\le_S$ is antisymmetric, we
  assume $F \le_S G$ and $G \le_S F$, witnessed by
  $\phi_1 : F \to G$ and $\phi_2 : G \to F$. Then
  defining $\phi := \phi_2 \circ \phi_1$,
   we obtain
   an injective map $\phi : F \to F$ such that $f \le \phi (f)$
   for all $f \in F$. We claim that
  $\phi (f) = f$ for all $f \in F$. Let $f$ be minimal
  in $F$ with respect to $\le_P$ such that $f \neq \phi(f)$.
  Then $f <_P \phi (f)$. Since $F$ is finite, $\phi$ is
  surjective, and thus there is $g \in F$ with
  $\phi (g) = f$. Since $\phi (f) >_P f$, we then
  have $g \neq f$ and therefore since $g \le_P \phi (g)$,
  we have $g <_P f$. Since we also have  $g \neq \phi (g)$,
  the polynomial vector $g$ contradicts the minimality of $f$. Therefore,
  $\phi$ is the identity map on $F$.
  Hence from $\phi_2 \circ \phi_1 = \operatorname{id}$, we obtain
  that for each $f \in F$, we have 
  $f \le_P \phi_1 (f) \le_P \phi_2 (\phi_1 (f)) = f$, which
  implies $\phi_1 (f) = f$ for all $f \in F$. Thus $\phi_1$
  is the identity mapping, which implies $F \subseteq G$.
  Since $F$ and $G$ have the same number of elements, this
  implies $F = G$, completing the proof that $\le_S$ is
  antisymmetric.

  Now let $(F_i)_{i \in \N}$ be such that
  $F_i >_S F_{i+1}$ for all $i \in \N$, and we choose such
  a chain for which $\# F_1$ is minimal.
  By this minimality, we then
  have $\# F_1 = \# F_i$ for all $i \in \N$.
  Let $\phi_i$ be an injective mapping
  from $F_{i+1}$ to $F_i$ with $f \le \phi_i (f)$ for all
  $f \in F_{i+1}$.
  Because of $\# F_i = \# F_{i+1}$, the mapping $\phi_i$ is
  bijective, and we have $\phi_i^{-1} (\phi_i(x)) = x \le_P \phi_i (x)$
  for all $x \in F_{i+1}$, and thus
  $\phi_i^{-1} (y) \le_P y$ for all $y \in F_i$.
  Let $\psi_i := \phi_{i}^{-1} \circ \cdots \circ \phi_2^{-1} \circ \phi_1^{-1}$,
  and fix $g \in F_1$. Then $(\psi_i (g))_{i \in \N_0}$ is
  a decreasing sequence in $(\rvxk, \le_P)$, and therefore,
  there is $n_1 \in \N$ such that
  $(\psi_i (g))_{i \in \N}$ is constant.
  Let $G_i := F_i \setminus \{\psi_{i-1} (g)\}$. The mappings
  $\phi_i \setminus \{ (\psi_i (g), \psi_{i-1} (g)) \}$ witness
  that $G_{i+1} \le_S G_i$ for all $i \in \N$.
  Hence by the minimality of $\# F_1$, the sequence $(G_i)_{i \in \N}$
  is constant from some $n_2$ onwards.
  Hence from $\max (n_1, n_2)$ onwards,
  $(F_i)_{i \in \N}$ is constant, a contradiction. 
\end{proof}

In proving that $S$-polynomial vectors that have a strong representation
still have a strong representation after applying
\textsc{SoftlyReduce} or \textsc{Normalize}, we have needed
the following lemma:
\begin{lem} \label{lem:transssr}
  Let $R$ be a Euclidean domain, let $\le$ be an admissible
  term order of $R$, and let
  $F, G, H \subseteq \rxnk$.
  We assume that every $f \in F$ has a strong standard representation
  by $G$ and that  every $g \in G$ has a strong standard representation
  by $H$. Then every $f \in F$ has a strong standard representation
  by~$H$.
\end{lem}
\begin{proof}
  If $f = \sum_{i=1}^N a_i n_i g_i$ is a strong standard representation
  of $f$ by $G$ and $g_i = \sum_{j=1}^{M_i} b_{i,j} m_{i,j} h_{i,j}$ is
  a strong standard representation of $g_i$ by $H$, then
  $f = \sum_{i=1} \sum_{j=1}^{M_i} a_i b_{i,j} \, (n_i m_{i,j}) \, h_{i,j}$
  is a representation of $f$ by $H$. To show that it is a strong
  standard representation, we observe that 
  $\md (n_1 m_{1,1} h_{1,1}) = \md (n_1) + \md (m_{1,1} h_{1,1}) = \md (n_1) + \md (g_1)$ because $g_1 = \sum_{j=1}^{M_1} b_{1,j} m_{1,j} h_{1,j}$ is
  a strong standard representation of $g_1$ by $H$.
  Furthermore, we have
  $\md (n_1) + \md (g_1) = \md (n_1g_1) = \md (f)$ because of the standard representation
  of $f$.
  Similarly, we see that for $(i,j) \neq (1,1)$, we have
  $\md (n_i m_{i,j} h_{i,j}) < \md (f)$.
\end{proof}

\section*{Acknowledgement}
I have discussed these topics with a number of colleagues at
our university. M.~Kauers provided several useful comments on the
manuscript. I thank these persons for their help.
\bibliography{gz18}
\end{document}